\documentclass[12pt]{article}

\usepackage{graphicx,tabularx,amsmath,mathtools,amsfonts,bm,amssymb,setspace,color,float,caption,subcaption,verbatim}
\usepackage[pdfpagemode=UseNone,pdfstartview={XYZ null null null}]{hyperref}
\usepackage[longnamesfirst,numbers,square]{natbib}
\usepackage[margin=0.5in]{geometry} 
\usepackage{ulem}
\usepackage{epsfig, graphicx}
\usepackage{amsmath}

\usepackage{latexsym}
\usepackage{amstext}
\usepackage{epsfig, graphicx}
\usepackage{latexsym}
\usepackage{amstext}
\usepackage{amssymb}
\usepackage{graphics}
\usepackage{color}

\newcommand{\p}{\partial}

\newcommand{\beq}{\begin{eqnarray}}
\newcommand{\beqq}{\begin{eqnarray*}}
\newcommand{\eeq}{\end{eqnarray}}
\newcommand{\eeqq}{\end{eqnarray*}}
\newcommand{\eps}{\varepsilon}

\newcommand{\R}{\mathbb{R}}

\newcommand{\ubar}{\overline{u}}

\definecolor{red}{rgb}{1,0,0}

\newtheorem{pres}{Principal Result}
\newtheorem{res}{Result}

\begin{document}
\title{Modeling and asymptotic analysis of the concentration difference in a nanoregion between an influx and outflux diffusion across narrow windows}
\author{{F. Paquin-Lefebvre$^{1}$ and D. Holcman$^{1,2}$} \footnote{$^{1}$ Applied Mathematics and Computational Biology, IBENS, Ecole Normale Sup\'erieure, 75005 Paris, France. $^2$ DAMTP, University of Cambridge, DAMTP and Churchill College CB30DS, United Kingdom.}}
\date{\today}
\maketitle
\begin{abstract}
When a flux of Brownian particles is injected in a narrow window located on the surface of a bounded domain, these particles diffuse and can eventually escape through a cluster of narrow windows. At steady-state, we compute asymptotically the distribution of concentration between the different windows. The solution is obtained by solving Laplace's equation using Green's function techniques and second order asymptotic analysis, and depends on the influx amplitude, the diffusion properties as well as the geometrical organization of all the windows, such as their distances and the mean curvature. We explore the range of validity of the present asymptotic expansions using numerical simulations of the mixed boundary value problem. Finally, we introduce a length scale to estimate how deep inside a domain a local diffusion current can spread. We discuss some applications in biophysics.
\end{abstract}


\section{Introduction}
How far inside a domain an influx of Brownian particles entering through a narrow window can perturb the steady-state bulk when the particles can escape through a neighboring window? We study here the properties of diffusion inside a bounded domain between two narrow windows. In the classical narrow escape theory \cite{ward93,JSP2004,HolcmanPNAS2007,pillay2010,cheviakov2010,HolcmanSchuss2015}, a stochastic particle initially distributed at a point or uniformly inside a bounded domain escapes through one of several narrow windows located on the surface. For that problem, asymptotic analysis and numerical simulations allow to study the relative contribution of the geometrical parameters on the mean escape time. Recently, the question of escape time has been extended to the fastest particles among many leading to an asymptotic formula that depends on the distance between the initial position and the target and the reciprocal of the logarithm of the number of particles \cite{Basna2018JNL,JCP2020,PLA2019}. Interestingly, when the distribution of initial particle overlays with the absorbing window, much faster escape times are expected as described in \cite{Weiss,Suney2021}.\\
In the present article, we do not study the escape when a particle is placed inside the domain, but consider a flow of particles entering through one window that spread inside a domain such as a ball, and can leave through another ensemble of target windows (Fig.~\ref{fig:fig1}A-C). This situation is inspired by the context of cell physiology. Indeed, the membrane potential is modulated due to the influx and efflux of ions across neighboring channels \cite{Hille,eisenberg,bezanilla2000,bezanilla2008}. These local influxes modify the local ionic concentrations and the electric field in small nanoregions, difficult to access experimentally. We focus on estimating the size of these regions where the concentration is perturbed, in the diffusion approximation for the motion of charged ions. In particular, to evaluate the size of this nanodomain, we introduce a new length measured by the flow line going from the center of a window where particles are injected to the center of the neighboring window where the flow is expelled.\\
We summarize the main asymptotic formulas that we derived for the difference of concentrations $c(y_1) - c(y_2)$ between two windows centered at points $y_1$ and $y_2$, with $c(y)$ solution of the Laplace's equation (see eq.~\eqref{eq:sspde}) when an influx $I$ is entering at $y_1$. Other parameters are the distance $L=\|y_1-y_2\|$ between the center of two circular windows and the common radius $\epsilon$ for all windows. Our main result is for a Brownian motion with diffusion coefficient $D$:
\begin{figure}[H]
\centering
\includegraphics[width=\linewidth]{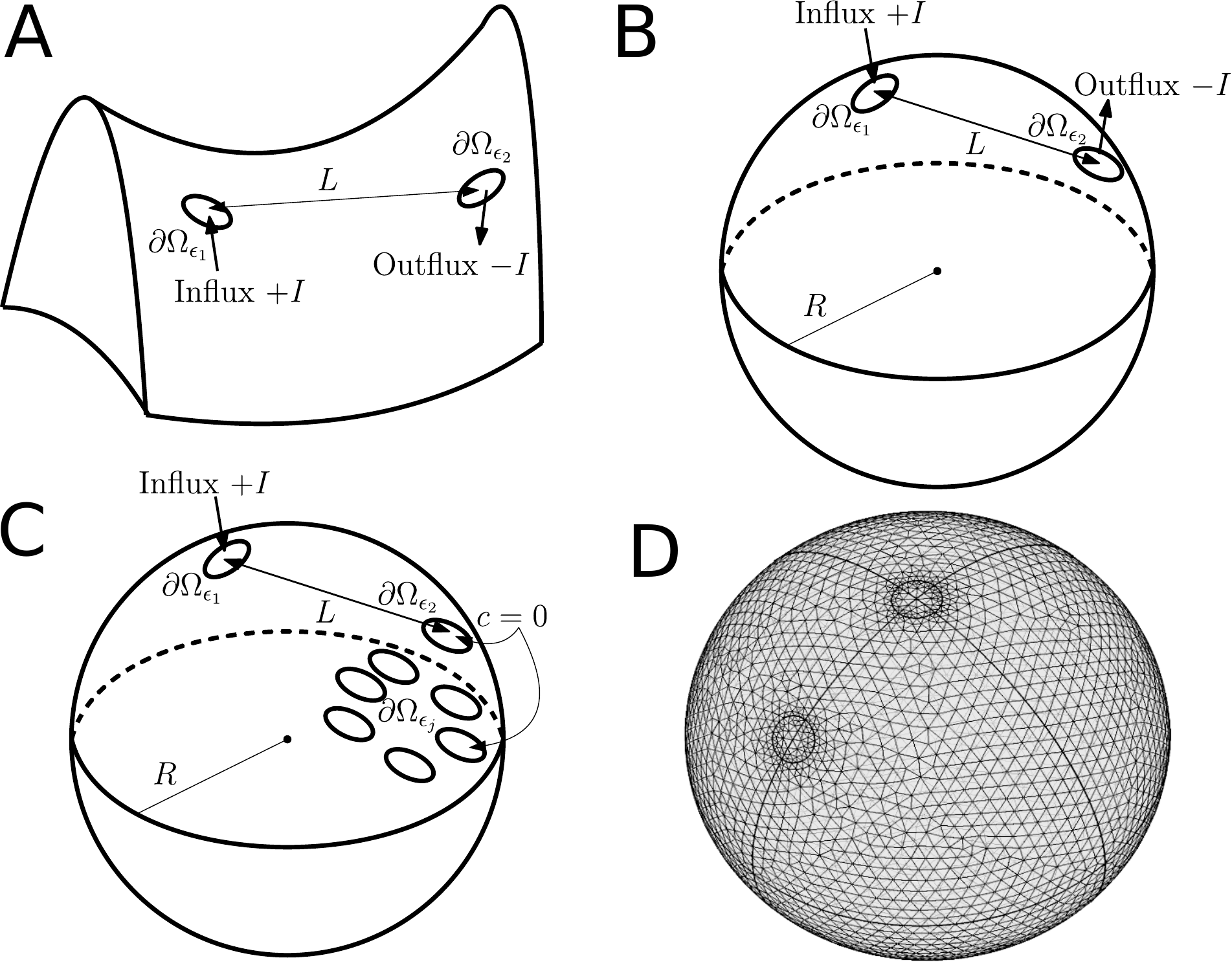}
\caption{{\bf Schematic representation of the influx-outflux} {\bf A:} Influx/outflux diffusion on an arbitrary domain. {\bf B:} Spherical geometry with Neumann boundary conditions. {\bf C:} Spherical geometry with mixed Neumann-Dirichlet boundary conditions. {\bf D:} Illustration of the a free tetrahedral meshing of a unit ball with two windows generated by COMSOL Multiphysics Version 5.2a \cite{comsol} used to solve numerically the distribution of the concentration.}
\label{fig:fig1}
\end{figure}
\paragraph{Two-term asymptotic formulas for an arbitrary closed and bounded domain $\Omega$ with smooth boundary $\p\Omega$.} When the isoperimetric ratio satisfies $|\p \Omega|/|\Omega|^{2/3} \sim O(1)$, this excludes long corridors with narrow necks of the order of $\epsilon$. We estimate the difference of concentration between two windows located on the boundary $\p\Omega$ which depends on the geometrical characteristic length-scale of the domain $R \approx \sqrt{|\p\Omega|}$, when the radii $R(P)$ of osculating spheres at all points $P$ on the boundary are of the same order as $R$ ($kR\leq \min_{P \in \p\Omega} R(P)\leq \max_{P \in \p \Omega} R(P) \leq K R$ for two order one constants $k$ and $K$). When $\frac{\epsilon}{R} \ll 1$, we have:
\begin{itemize}
\item for two small windows of radius $\epsilon$ with Neumann boundary conditions $(+I)$ at window $\p \Omega_{\epsilon_1}$ (influx) and $(-I)$ at window $\p \Omega_{\epsilon_2}$ (outflux), then
\beq
c(y_1) - c(y_2) = \frac{\epsilon I}{D}\left(2 - \frac{H\left(y_1\right) + H\left(y_2\right)}{4}{\epsilon}\log\left(\frac{\epsilon}{R}\right) + O\left(\frac{\epsilon}{R}\right) \right)\,,
\eeq
where $H(x)$ is the mean curvature at a point $x$.
\item for two small windows of radius with an influx condition $(+I)$ on $\p \Omega_{\epsilon_1}$ and absorbing boundary condition on $\p \Omega_{\epsilon_2}$, then
\beq
c(y_1) = \frac{\epsilon I}{D}\left(1 + \frac{\pi}{4} - \frac{H\left(y_1\right) + H\left(y_2\right)}{4}{\epsilon}\log\left(\frac{\epsilon}{R}\right) + O\left(\frac{\epsilon}{R}\right) \right)\,,
\eeq
\item for an influx condition on $\p \Omega_{\epsilon_1}$ surrounded by $N-1$ circular patches $\p \Omega_{\epsilon_j}$ for $j=2,\ldots,N$ with absorbing boundary conditions, then
\beq
c(y_1) = \frac{\epsilon I}{D}\left(1 + \frac{\pi}{4(N-1)} - \left(H\left(y_1\right) + \frac{1}{(N-1)^2}\sum_{i=2}^N H\left(y_i\right) \right) \frac{\epsilon}{4}\log\left(\frac{\epsilon}{R}\right) + O\left(\frac{\epsilon}{R}\right) \right)\,.
\eeq
\end{itemize}
For a ball $B(R)$ of radius $R$:
\begin{itemize}
\item for two windows with a flux boundary conditions (influx $+I$ in $\p \Omega_{\epsilon_1}$ and outflux $-I$ in $\p \Omega_{\epsilon_2}$),
\beq
c(y_1) - c(y_2) = \frac{\epsilon I}{D}\left(2 - \frac{\epsilon}{2R}\log\left(\frac{\epsilon}{R}\right) + \left( \frac{1}{4} - \frac{R}{L} + \frac{1}{2}\log\left(\frac{L^2}{2R^2} + \frac{L}{R}\right)\right)\frac{\epsilon}{R} + O\left(\left(\frac{\epsilon}{R}\right)^2\right) \right)\,.
\eeq
\item when the influx condition is prescribed on $\p \Omega_{\epsilon_1}$, with an absorbing boundary condition on $\p \Omega_{\epsilon_2}$,
\beq
c(y_1) = \frac{\epsilon I}{D}\left(1+\frac{\pi}{4} - \frac{\epsilon}{2R}\log\left(\frac{\epsilon}{R}\right) + \left( \frac{3}{8} - \frac{\log(2)}{4} - \frac{R}{L} + \frac{1}{2}\log\left(\frac{L^2}{2R^2} + \frac{L}{R}\right) \right)\frac{\epsilon}{R} + O\left(\left(\frac{\epsilon}{R}\right)^2\right) \right)\,.
\eeq
\item when an influx condition is on $\p \Omega_{\epsilon_1}$ and there are $N-1$ circular windows with absorbing boundary conditions,
\beq
\begin{split}
& c(y_1) = \frac{\epsilon I}{(N-1)D}\left(N-1 + \frac{\pi}{4} - \frac{N\epsilon}{4R}\log\left(\frac{\epsilon}{R}\right) + \left( \frac{N+1}{8} - \frac{\log(2)}{4} \right. \right. \\
& \left. \left. + \frac{1}{N-1}\sum_{i=2}^N\sum_{j=i+1}^N \left( \frac{R}{L_{ij}} - \frac{1}{2}\log\left( \frac{L_{ij}^2}{2R^2} + \frac{L_{ij}}{R} \right) \right) - \sum_{j=2}^N \left( \frac{R}{L_{1j}} - \frac{1}{2}\log\left( \frac{L_{1j}^2}{2R^2} + \frac{L_{1j}}{R} \right) \right)\right)\frac{\epsilon}{R} \right. \\
& \left. + O\left(\frac{\epsilon^2}{R^2}\right)\right)\,,
\end{split}
\eeq
where $L_{ij}= \|y_i - y_j\|$.
\end{itemize}
The manuscript is organized as follows: in Part 1, we introduce the diffusion model equations. In Part 2, we derive the asymptotic formula for the difference of concentration between two neighboring windows in two cases: for a Dirichlet boundary condition on the exiting hole or Neumann condition with a flux which is the same with opposite sign as the influx. Numerical solutions are performed by COMSOL \cite{comsol} using a mesh decomposition as shown in Fig.\ref{fig:fig1}D. In Part 3, we discuss the case of the several windows and finally in Part 4, we define the length of penetration using the flow line.
\section{Modeling diffusion from the fluxes to small circular patches on a general domain}
We consider a bounded domain $\Omega$ with boundary $\p\Omega$ divided into $N-1$ small and disjoint absorbing circular windows $\p\Omega_{\epsilon_j}\ (j=2,\ldots,N)$, centered at $y_j$ and of radius $\epsilon$. Each has area $|\p \Omega_{\epsilon_j}|=O(\epsilon^2)$ (in three dimensions). The total absorbing boundary is
\beq
\p\Omega_a= \cup_{j=2}^N \p \Omega_{\epsilon_j}.
\eeq
On the boundary there is also an additional window $\p\Omega_{\epsilon_1}$ that receives an influx of Brownian particles with a steady-state amplitude current $I$. The windows are not necessarily far apart so that non-linear effects~\citep{holcman2008} could be expected. The remaining boundary surface, $\p\Omega_r=\p\Omega\setminus\{\p\Omega_{\epsilon_1}\cup\p\Omega_a\}$, is reflective for the diffusing particles. Other models are possible, and instead of absorbing boundary conditions we could consider partially absorbing (Robin) boundary conditions \citep{grebenkov2020single}. The concentration $c(y,t)$ of Brownian particles at position $y$ at time $t$ satisfies the diffusion equation
\beq\label{IBVP_mixedND}
\begin{split}
\frac{\p c(y,t)}{\p t} &= D \Delta c(y,t)\quad\hbox{\rm for } y \in \Omega\,,\\
\frac{\p c(y,t)}{\p n} &= 0\quad \hbox{\rm for } y \in\p\Omega_r\,,\\
D\frac{\p c(y,t)}{\p n} &= I\quad \hbox{\rm for } y \in\p\Omega_{\epsilon_1}\,,\\
c(y,t) &= 0\quad \hbox{\rm for } y \in \p\Omega_a\,,
\end{split}
\eeq
where $D$ is the diffusion coefficient and $n$ is the outward normal to the boundary. When there are $N=2$ windows, with $\p \Omega_{\eps_1}$ receiving an influx $I$ and $\p\Omega_{\epsilon_2}$ emitting an outflux of opposite amplitude $-I$, the diffusion model becomes
\beq\label{eq:IBVP_nbc}
\begin{split}
\frac{\p c(y,t)}{\p t} &= D \Delta c(y,t)\quad\hbox{\rm for } y \in \Omega\,,\\
\frac{\p c(y,t)}{\p n} &= 0\quad \hbox{\rm for } y \in\p\Omega_r\,,\\
D\frac{\p c(y,t)}{\p n} &= I\quad \hbox{\rm for } y \in\p\Omega_{\epsilon_1}\,,\\
D\frac{\p c(y,t)}{\p n} &= -I\quad \hbox{\rm for } y \in\p\Omega_{\epsilon_2}.
\end{split}
\eeq
The domain $\Omega$ could contain an initial distribution of Brownian particles at equilibrium, thus if there are initially $Q$ particles, the density will be $\rho_0(y)=\frac{Q}{V}$ with volume $V=\text{vol}(\Omega)$. Thus in that case, we will be interested in a steady-state perturbation of the concentration $\rho(y)=\rho_0(y)+\rho_1(y)+\ldots$, where $\rho_1(y)$ is the first order term due to local changes in the concentration induced by the influx-outfluxes.\\
At steady-state we have $\frac{\p c(y,t)}{\p t}=0$, and thus for two windows receiving and emitting a flux of particles of amplitude $I$, the steady-state concentration $c$ satisfies the Laplace's equation
\beq\label{eq:sspde}
\begin{split}
D\Delta c &= 0\,, \quad y \in \Omega\, \\
D\frac{\p  c}{\p  n} &= 0\,, \quad y \in \p \Omega_r \,, \quad D\frac{\p  c}{\p  n} = I\,, \quad y \in \p\Omega_{\epsilon_1}\,, \quad D\frac{\p  c}{\p  n} = -I\,, \quad y \in \p\Omega_{\epsilon_2}\,, \\
\end{split}
\eeq
The solution is defined up to an additive constant, which will cancel out by taking the difference of concentration between two points. The non-dimensionalization variable is given by
\beq
u(x)=\frac{c(xR)D}{IR}\,, \quad \text{with} \quad x = \frac{y}{R} \in \tilde{\Omega} \,.
\eeq
Thus $u(x)$ is solution of
\beq\label{eq:pdenon}
\begin{split}
\Delta u &= 0\,, \quad x \in \Omega\,, \\
\frac{\p  u}{\p  n} &= 0\,, \quad x \in \p \Omega_r \,, \quad \frac{\p  u}{\p  n} = 1\,, \quad x \in \p\Omega_{\eps_1}\,, \quad \frac{\p  u}{\p  n} = -1\,, \quad x \in \p\Omega_{\eps_2}\,.
\end{split}
\eeq
where we dropped the tilde symbol above $\Omega$ to simplify notations. The remaining two parameters consist of the normalized window radius $\eps \equiv \epsilon/R$, in addition of the normalized distance between the centers of the windows $l \equiv L/R = \|x_1-x_2\|$.
\subsection{Asymptotic solution using Neumann-Green function}\label{sec:nbc}
To obtain an explicit solution $u(x)$ of eq.~\eqref{eq:pdenon} we use any of the Neumann-Green's function $G_s(x;y)$ which satisfies
\beq\label{eq:greenfun}
\Delta G_s(x;y) = \frac{1}{|\Omega|}\,, \quad x \in \Omega\,; \quad \frac{\p  G_s}{\p  n} = \delta(x-y)\,, \quad x \in \p  \Omega\,; 
\eeq
with $y\in \p  \Omega$. For arbitrary domains the Neumann's function has the following expansion near the singular diagonal $x = y$ \cite{pre2008,HolcmanSchuss2015,gomez2015}
\beq\label{eq:gengf}
G_s(x;y) = \frac{1}{2\pi \|x-y\|} - \frac{1}{4\pi}H(y)\log\left(\|x-y\|\right) + V_s(x;y)\,,
\eeq
where $H(y)$ is the mean curvature computed at $y \in \p \Omega$, and $V_s(x;y)$ is a bounded and regular function of $x$ and $y$ in $\Omega$. Using Green's second identity we have
\beq\label{eq:greenid}
\int_{\Omega} \left( G_s\Delta u - u \Delta G_s \right)dx = \int_{\p \Omega} \left( G_s\frac{\p u}{\p n} - u \frac{\p G_s}{\p n} \right)dx\,,
\eeq
and substituting \eqref{eq:pdenon} and \eqref{eq:greenfun} within \eqref{eq:greenid}, we obtain for two windows
\beq\label{eq:ubdy}
u(y) = \ubar + \int_{\p \Omega_{\eps_1}}G_s(x;y) \frac{\p u}{\p n} dx + \int_{\p \Omega_{\eps_2}}G_s(x;y) \frac{\p u}{\p n}dx = \ubar + \int_{\p \Omega_{\eps_1}}G_s(x;y)dx - \int_{\p \Omega_{\eps_2}}G_s(x;y)dx\,,
\eeq
where
\beq
\ubar \equiv \frac{1}{|\Omega|} \int_\Omega u dx\,.
\eeq
We shall next compute each integral term by considering $y=x_1$ and $y=x_2$ in \eqref{eq:ubdy}. Taking the difference, we get
\beq\label{eq:soldiffgen}
u(x_1) - u(x_2) = \int_{\p  \Omega_{\eps_1}}G_s(x;x_1)dx + \int_{\p  \Omega_{\eps_2}}G_s(x;x_2)dx - \int_{\p  \Omega_{\eps_2}}G_s(x;x_1)dx - \int_{\p  \Omega_{\eps_1}}G_s(x;x_2)dx \,.
\eeq
To estimate the integrals of the Neumann Green's function near a singularity, we use \eqref{eq:gengf} to approximate $G_s(x;x_j)$ for $x \in \p \Omega_{\eps_j}$ as
\beq\label{eq:gdec}
G_s(x;x_j) \sim g_j(r) + v_j  \,, \quad r = \|x-x_j\| < \eps \,, \quad j = 1,\,2,
\eeq
where $v_j \equiv V_s(x_j;x_j)$ is constant and $g_j(r)$ is the sum of the Coulomb and logarithmic singularities,
\beq
g_j(r) \equiv \frac{1}{2\pi r} - \frac{1}{4\pi}H(x_j)\log\left(r\right)\,.
\eeq
A direct integration using polar coordinates yields
\begin{align}
\int_{\p \Omega_{\eps_j}}G_s(x;x_j)dx &= 2\pi\int_0^\eps \left(g_j(r) + v_j \right)rdr = \eps - \frac{H(x_j)}{2}\eps^2\int_0^1 \log(\eps u) u du + v_j \pi \eps^2\,, \nonumber \\
&= \eps - \frac{H(x_j)}{4} \eps^2\log(\eps) + \frac{H(x_j)}{8}\eps^2 + v_j \pi \eps^2\,. \label{eq:intref}
\end{align}
The last two terms in \eqref{eq:soldiffgen} correspond to integrals of the Neumann's function away from the singularity, and can be approximated as
\beq
\int_{\p  \Omega_{\eps_i}} G_s(x;x_j)dx \sim G_s(x_i;x_j)\pi\eps^2\,, \quad i \neq j\,,
\eeq
Finally, by adding up the different terms in \eqref{eq:soldiffgen}, we obtain that the concentration difference between the center of each window is given by
\beq\label{eq:NeumannSolDiff}
u(x_1) - u(x_2) = 2\eps - \frac{H(x_1) + H(x_2)}{4}\eps^2\log(\eps) + \left( \frac{H(x_1)+H(x_2)}{8} + \pi\left(v_1 + v_2 - 2G_s(x_1;x_2)\right)\right)\eps^2.
\eeq
To conclude, formula \eqref{eq:NeumannSolDiff} is an asymptotic approximation for the concentration difference between the flux receiving and emitting windows. The circular patches of radius $\eps$ could at most be tangent on a smooth boundary. The second term originates from the logarithmic singularity of the 3D Neumann Green's function. Finally, the first two terms are independent of the distance between each window and do not require an explicit solution for the Green's function. Therefore, dropping the higher-order quadratic term yields a general two-term asymptotic formula for the difference of concentration that only requires the knowledge of the mean curvature at the center of each hole. In summary,
\begin{res}
For two circular windows of small radius $\eps \ll 1$ and separated by an order one distance $\|x_1 - x_2\| \sim O(1)$, the normalized solution $u(x)$ of \eqref{eq:pdenon} in the domain $\Omega$ has a two-term asymptotic approximation for the concentration drop between the two windows given by
\beq\label{eq:diffugen}
u(x_1) - u(x_2) = 2\eps - \frac{H(x_1) + H(x_2)}{4}\eps^2\log(\eps) + O(\eps^2)\,,
\eeq
where $H(x_1)$ and $H(x_2)$ are the mean curvature at $x_1$ and $x_2$ on the boundary $\p\Omega$.
\end{res}
\subsection{Difference of concentration with a point inside the domain}
To calculate the difference of concentration with a point inside the domain, we use a different Neumann's Green function \cite{JSP2004,HolcmanSchuss2015}, defined as
\beq
\label{eq:diff_G}
\begin{aligned}
\Delta G(x;y) = -\delta(x-y) \quad & x \in \Omega,\\
\frac{\p  G (x,y)}{\p  n} = -\frac{1}{|\p \Omega|} \quad & x \in \p  \Omega\,,
\end{aligned}
\eeq
where the singularity $y$ is not anymore restricted to the boundary. Using Green's second identity, we obtain for $y \in \Omega$,
\beq
u(y) = \frac{1}{|\p\Omega|} \int_{\p \Omega} u dx + \int_{\p \Omega_{\eps_1}} G(x;y) dx - \int_{\p \Omega_{\eps_2}} G(x;y) dx\,.
\eeq
We get
\beq
u(x_1)- u(y) = \int_{\p  \Omega_{\eps_1}} G(x;x_1) dx-\int_{\p \Omega_{\eps_2}} G(x;x_1) dx - \int_{\p \Omega_{\eps_1}} G(x;y) dx+\int_{\p \Omega_{\eps_2}} G(x;y) dx\,,
\eeq
which becomes,
\beq
u(x_1)- u(y) = \eps - \frac{H(x_1)}{4} \eps^2\log(\eps) + \frac{H(x_1)}{8}\eps^2 + v_1 \pi \eps^2 - G(x_2;x_1) |\p \Omega_{\eps_1}| - G(x_1;y) |\p \Omega_{\eps_1}|  + G(x_2;y)|\p \Omega_{\eps_2}|\,.
\eeq
Thus when the point $y$ is located at an intermediate position between $x_1$ and $x_2$, such that $G(x_1;x_2)\ll G(y;x_2), G(y;x_1)$, we can use the expansion of the Neumann-Green's function,
\beq
G(y;x_j) = \frac{1}{2\pi \|y-x_j\|} - \frac{1}{4\pi}H(x_j)\log\left(\|y-x_j\|\right) + V_s(y;x_j)\,, \quad j = 1,\,2\,,
\eeq
to obtain the following asymptotic expansion for the difference of concentration,
\begin{align}
u(x_1) - u(y) &= \eps - \frac{H(x_1)}{4} \eps^2\log(\eps) + \left(\frac{H(x_1)}{8} + v_1\pi - \pi G(x_2;x_1) - \frac{1}{2\|y-x_1\|} + \frac{1}{2\|y-x_2\|} \right. \nonumber \\
& \left. + \frac{1}{4\pi}H(x_1)\log\left(\|y-x_1\|\right) - \frac{1}{4\pi}H(x_2)\log\left(\|y-x_2\|\right) - V_s(y;x_1) + V_s(y;x_2) \right) \eps^2\,.
\end{align}
\subsection{Difference of concentration when the two windows are not far apart}\label{sec:closeholes}
We compute here the difference of concentration when the distance between the two windows $\p \Omega_{\eps_1}$ and $\p \Omega_{\eps_2}$ is small, of the order $l = \|x_1-x_2\| \sim O(\eps)$. Even if the windows are not necessarily in the same tangent plane on a Riemannian  surface, they are located nearby and thus their projections deviate at second order. The projected two windows are located on the same tangent plane. The calculation starts from expression \eqref{eq:soldiffgen}, and the first two terms are calculated exactly as before. However for the last two terms we use the law of cosines to express the distance from the singularity, which gives
\beq
u(x_1) - u(x_2) = 2\pi\sum_{j=1}^2 \int_0^\eps g_j(r)rdr - \sum_{j=1}^2\int_0^{2\pi}\int_0^\eps g_j(\sqrt{l^2 + r^2 -2lr\cos(\theta)})rdrd\theta\,,
\eeq
where the regular parts of the Green's function have been omitted since they disappear in the difference. After a change of variable, we get
\beq\label{eq:lsmall1}
u(x_1) - u(x_2) = 2\pi \eps^2 \sum_{j=1}^2 \int_0^\eps g_j(u)udu - \eps^2\sum_{j=1}^2\int_0^{2\pi}\int_0^\eps g_j(\sqrt{l^2 + (u\eps)^2 -2lu\eps\cos(\theta)})udud\theta\,.
\eeq
Next, using our result in \eqref{eq:intref}, the first sum evaluates as
\beq
2\pi\eps^2\sum_{j=1}^2 \int_0^\eps g_j(u)udu = 2\eps - \frac{H(x_1)+H(x_2)}{4} \eps^2\log(\eps) + \frac{H(x_1)+H(x_2)}{8}\eps^2\,.
\eeq
Using $l = \eta \eps$, with $\eta \geq 2$ to avoid overlapping, and defining $\omega(u,\theta,\eta)$ as
\beq
\omega(u,\theta,\eta) = \eta^2 - 2\eta u\cos(\theta) + u^2\,,
\eeq
we obtain that
\beq
g_j(\sqrt{l^2 + r^2 -2lr\cos(\theta)}) = g_j(\eps\sqrt{\omega(u,\theta,\eta)})\,,
\eeq
and then calculating the integral yields
\begin{align*}
&\int_0^{2\pi} \int_0^1 g_j(\eps\sqrt{\omega(u,\theta,\eta)})udud\theta = \int_0^{2\pi} \int_0^1 \left( \frac{u}{2\pi \eps \sqrt{\omega(u,\theta,\eta)}} - \frac{H(x_j)}{4\pi}u\log\left(\eps\sqrt{\omega(u,\theta,\eta)}\right) \right)dud\theta\,, \\
&= \frac{1}{2\pi\eps} \int_0^{2\pi} \int_0^1 \frac{udud\theta}{\sqrt{\omega(u,\theta,\eta)}} - \frac{H(x_j)}{4\pi}\log(\eps) - \frac{H(x_j)}{8\pi}\int_0^{2\pi}\int_0^1 \log\left(\omega(u,\theta,\eta)\right)udud\theta\,.
\end{align*}
Then, by switching the order of integration,
\beq
\frac{1}{2\pi\eps}\int_0^{2\pi}\int_0^1\frac{udu d\theta}{\sqrt{\omega(u,\theta,\eta)}} = \frac{1}{2\pi\eps}\int_0^1\int_0^{2\pi}\frac{u d\theta du}{\sqrt{\eta^2 - 2\eta u\cos(\theta) + u^2}} = \frac{2}{\pi\eps}\int_0^1 \frac{u}{\eta + u}K\left(\frac{2\sqrt{\eta u}}{\eta + u}\right)du\,,
\eeq
where $K(\cdot)$ is the complete elliptic integral of the first kind. Finally, by summing all the terms in \eqref{eq:lsmall1}, the $O(\eps^2\log(\eps))$ terms cancel and we obtain the following expression
\begin{align}\label{eq:lsmall2}
&u(x_1) - u(x_2) = \nonumber \\
& \left(2 - \frac{4}{\pi}\int_0^1 \frac{u}{\eta + u}K\left(\frac{2\sqrt{\eta u}}{\eta + u}\right)du\right)\eps + \left(\frac{1}{8} + \frac{1}{8\pi}\int_0^{2\pi}\int_0^1 \log\left(\omega(u,\theta,\eta)\right)udud\theta\right)\left(H(x_1) + H(x_2)\right)\eps^2 \,,
\end{align}
which gives, at leading-order,
\begin{align}\label{eq:lsmall3}
u(x_1) - u(x_2) = \left(2 - \frac{4}{\pi}\int_0^1 \frac{u}{\eta + u}K\left(\frac{2\sqrt{\eta u}}{\eta + u}\right)du\right)\eps + O(\eps^2) \,.
\end{align}
when $\eta = 2$ in \eqref{eq:lsmall2}, the two windows are tangent holes and we obtain a correction to the leading-order term for the concentration difference: by numerically integrating \eqref{eq:lsmall3} using the built-in quadrature routine from Matlab, we find
\begin{align}\label{eq:lsmallnum}
\left(u(x_1) - u(x_2)\right)|_{\eta = 2} \approx 1.41676\eps + O(\eps^2) < 2\eps + O(\eps^2) \,.
\end{align}
Thus, at leading-order order we obtain weaker concentration differences than what the asymptotic theory predicts for well-separated windows. Alternatively when $\eta \gg 2$ the two windows are far apart, and we show below that we recover formula \eqref{eq:diffugen} as a limiting case. By using the expansion $K(\cdot) \approx \frac{\pi}{2}$ for small arguments of the elliptic integral, we obtain that
\beq
\frac{4}{\pi}\int_0^1 \frac{u}{\eta + u}K\left(\frac{2\sqrt{\eta u}}{\eta + u}\right)du \approx \frac{4}{\pi}\frac{\pi}{2}\int_0^1 \frac{udu}{\eta} = \frac{1}{\eta}\,.
\eeq
Similarly, if we express $\omega(u,\theta,\eta)$ as
\beq
\omega(u,\theta,\eta) = \eta^2 \left( 1 + \left(\frac{u}{\eta}\right)^2 - 2\frac{u}{\eta}\cos(\theta)\right)\,,
\eeq
then we can approximate the integral of the $\log$ term in \eqref{eq:lsmall2} as
\beq
\frac{1}{8\pi}\int_0^{2\pi}\int_0^1 \log\left(\omega(u,\theta,\eta)\right)udud\theta \approx \frac{1}{8\pi}\int_0^{2\pi}\int_0^1 \log\left(\eta^2\right)udud\theta = \frac{\log(\eta)}{4}\,.
\eeq
Now, recalling that $\eta = \frac{l}{\eps}$ we find that
\beq
u(x_1) - u(x_2) \approx 2\eps - \frac{\eps^2}{l} + \left( \frac{1}{8} + \frac{1}{4}\log\left(\frac{l}{\eps}\right) \right)\left(H(x_1) + H(x_2)\right)\eps^2 \,,
\eeq
which becomes, after rearranging the terms and dropping the $O(\eps^2)$ terms,
\beq
u(x_1) - u(x_2) = 2\eps - \frac{H(x_1) + H(x_2)}{4}\eps^2\log(\eps) + O(\eps^2)\,.
\eeq
Hence, we have recovered formula \eqref{eq:diffugen}, which was derived for well-separated circular windows.

\subsection{Asymptotic expression for the case of mixed Dirichlet-Neumann boundary conditions}\label{sec:mixedND}
We now extend the previous analysis to an arbitrary domain $\Omega$ whose boundary is punctured with several circular windows $\p \Omega_{\eps_j}$, each of radius $\eps$ and centered at a point $x_j$ for $j=1,2,\ldots,N$. An influx current is applied on the first hole $j=1$, while in the remaining windows we impose an absorbing boundary condition. When $\eps \ll 1$ and the windows are well-separated, the distances $l_{ij} \equiv \|x_i - x_j\|$ between each center is of $O(1)$ compared to $\eps$. For this case, the steady-state system \eqref{eq:pdenon} is given by
\beq\label{eq:pdenonN}
\begin{split}
\Delta u &= 0\,, \quad x \in \Omega\,, \\
\frac{\p  u}{\p  n} &= 0\,, \quad x \in \p\Omega_r\,, \quad \frac{\p  u}{\p  n} = 1\,, \quad x \in \p\Omega_{\eps_1}\,, \quad u = 0\,, \quad x \in \p\Omega_{\eps_j}\,, \quad j=2,\ldots,N\,. \\
\end{split}
\eeq
Green's identity applied to \eqref{eq:greenfun} and \eqref{eq:pdenonN} yields for any $y$ on the boundary $\p \Omega$
\beq\label{eq:ubdyN}
u(y) = \ubar + \int_{\p \Omega_{\eps_1}}G_s(x;y) dx + \sum_{j=2}^N \int_{\p \Omega_{\eps_j}}G_s(x;y)\frac{\p u}{\p n} dx.
\eeq
Our aim is to calculate the solution at a point $x_1$, since we already have $u(x_j) = 0$ for $j \neq 1$. Before proceeding, we recall that the normal derivative of the flux at each exiting window $\p \Omega_{\eps_j}$ $j = 2,\ldots,N$ is given by the classical Weber solution \cite{crank1975,holcman2008}
\beq\label{eq:weber}
\frac{\p u(x)}{\p n} = \frac{C_j}{\sqrt{\eps^2 - \|x-x_j\|^2}}\,, \quad \text{for} \quad \|x-x_j\| < \eps \,, \quad j \neq 1\,,
\eeq
where $C_j$ is a constant and thus the outflux $\Phi_j$ across $\p \Omega_{\eps_j}$ is defined as
\beq
\Phi_j \equiv \int_{\p\Omega_{\eps_j}} \frac{\p u(x)}{\p n} dx = 2\pi C_j \int_0^\eps \frac{rdr}{\sqrt{\eps^2 - r^2}} = 2\pi C_j \eps\,, \quad j \neq 1\,.
\eeq
The contribution from all exiting fluxes must compensate the total influx $\Phi_1 = \pi \eps^2$. Therefore, the comptability condition (divergence theorem) to \eqref{eq:pdenonN} insure that the unknown constants $C_j$ are linked by
\beq\label{eq:compat}
\sum_{j=2}^N C_j = - \frac{\eps}{2}\,.
\eeq
Next we set $y = x_i$ for $i = 2,\,\ldots,N$ within \eqref{eq:ubdyN}, and since $u(x_i) = 0$ from the boundary conditions, we obtain
\beq\label{eq:uxj1}
0 = \ubar + \int_{\p \Omega_{\eps_1}}G_s(x;x_i) dx + \sum_{j=2}^N \int_{\p \Omega_{\eps_j}}G_s(x;x_i)\frac{\p u}{\p n} dx \,.
\eeq
Away from a singularity $x_i$ on a patch $\p \Omega_{\eps_j}$ with $j \neq i$, we approximate the Green's function by a constant owing to the fact that $\|x_i - x_j\| \sim O(1)$ and $\eps \ll 1$ (see also section \S \ref{sec:nbc}). Therefore, using the boundary condition $\frac{\p u}{\p n} = 1$ on $\p\Omega_{\eps_1}$, we obtain
\beq
\int_{\p \Omega_{\eps_1}}G_s(x;x_i)\frac{\p u}{\p n} dx \sim \pi\eps^2 G_s(x_1;x_i)\,, \quad \int_{\p \Omega_{\eps_j}}G_s(x;x_i)\frac{\p u}{\p n} dx \sim 2\pi\eps G_s(x_j;x_i) C_j\,,
\eeq
which reduces \eqref{eq:uxj1} to
\beq\label{eq:uxj2}
0 = \ubar + \pi \eps^2 G_s(x_1;x_i) + 2\pi C_i \int_0^\eps \frac{(g_i(r)+v_i)rdr}{\sqrt{\eps^2-r^2}} + 2\pi\eps\sum_{\substack{j=2\\ j \neq i}}^N C_j G_s(x_j;x_i)\,,
\eeq
where $g_i(r)$ and $v_i = V_s(x_i;x_i)$ are the singular and regular parts of the Green's function as given in \eqref{eq:gdec}. The integral in expression \eqref{eq:uxj2} is directly computed as
\begin{align}
2 \pi \int_0^\eps \frac{(g_i(r) + v_i) r dr}{\sqrt{\eps^2 - r^2}} &= \int_0^1 \left( \frac{1}{\sqrt{1-u^2}} - \frac{\eps u}{2\sqrt{1-u^2}}H(x_i)\log\left(\eps u\right) + 2\pi \eps v_i\frac{u}{\sqrt{1-u^2}} \right) du\,, \nonumber \\
&= \frac{\pi}{2} + \left(-\frac{H(x_i)}{2}\log(\eps) + \frac{1-\log(2)}{2}H(x_i) + 2\pi v_i \right)\eps \label{eq:intabs} \,.
\end{align}
Next, by defining the quantity $d_i$ as
\beq\label{eq:d_i}
d_i \equiv -\frac{H(x_i)}{2}\log(\eps) + \frac{1-\log(2)}{2}H(x_i) + 2\pi v_i\,,
\eeq
and combining \eqref{eq:uxj2} with the compatibility condition \eqref{eq:compat}, we get the following system of equations for the constants $C_i$ and the average concentration $\ubar$,
\begin{align}
& C_i\left(\frac{\pi}{2} + d_i\eps \right) + 2\pi\eps\sum_{\substack{j=2\\ j \neq i}}^N C_j G_s(x_j;x_i) = - \ubar - \pi \eps^2 G_s(x_1;x_i) \quad \text{for} \quad i=2,\ldots,N \label{eq:pert1} \,, \\
& \sum_{j=2}^N C_j = - \frac{\eps}{2}\,.
\end{align}
We then write \eqref{eq:pert1} under matrix form as follows\,,
\beq\label{eq:pert2}
\left(\frac{\pi}{2}I_{N-1} + \eps \bm{M} \right) \bm{C} = -\ubar {1}_{N-1} - \pi \eps^2\bm{b}\,,
\eeq
where $I_{N-1}$ is the identity matrix of size $N-1$ and ${1}_{N-1}$ is a column vector with each element equal to 1. The matrix $\bm{M}$ is further defined by
\beq
\bm{M} =
\begingroup
\renewcommand*{\arraystretch}{1.5}
\begin{pmatrix}
d_2 & 2\pi G_s(x_3;x_2) & \cdots & 2\pi G_s(x_N;x_2) \\
\vdots & \vdots & \ddots & \vdots \\
2\pi G_s(x_2;x_N) & 2\pi G_s(x_3;x_N) & \cdots & d_N \\
\end{pmatrix}\,,
\endgroup
\eeq
and the vectors $\bm{C}$ and $\bm{b}$ by
\beq
\bm{C} =
\begin{pmatrix}
C_2 \\
\vdots \\
C_N
\end{pmatrix}\,,
\quad
\bm{b} =
\begin{pmatrix}
G_s(x_1;x_2) \\
\vdots \\
G_s(x_1;x_N)
\end{pmatrix}\,.
\eeq
The matrix equation \eqref{eq:pert2} is defined as the sum of the identity with an $O(\eps)$ perturbation, thus it is invertible and we get
\beq\label{eq:pert3}
\bm{C} =-\frac{2}{\pi}\left(I_{N-1} + \frac{2\eps} {\pi} \bm{M} \right)^{-1}\left(\ubar {1}_{N-1} + \pi \eps^2\bm{b}\right) = -\frac{2}{\pi}\left(\sum_{n=0}^{\infty} \left(-\frac{2\eps} {\pi}\right)^n \bm{M}^n\right)\left(\ubar {1}_{N-1} + \pi \eps^2\bm{b}\right)\,,
\eeq
and then keeping the first two terms of the geometric series yields
\beq\label{eq:twoTerms}
\bm{C} = -\frac{2}{\pi}\left(I_{N-1} - \frac{2\eps}{\pi} \bm{M}\right)\left( \ubar {1}_{N-1} + \pi \eps^2\bm{b} \right) =  -\frac{2}{\pi}\left( \ubar \left( {1}_{N-1} - \frac{2\eps}{\pi}\bm{M}{1}_{N-1} \right) + \pi \eps^2 b \right) + O(\eps^3)\,.
\eeq
Summing over all the rows of \eqref{eq:twoTerms} and using the compatibility condition \eqref{eq:compat} on the left-hand side leads to an equation for the average concentration $\ubar$,
\beq
-\frac{\eps}{2} = -\frac{2}{\pi}\left( \ubar(N-1) - \ubar\frac{2\eps}{\pi}\left( \sum_{i=2}^N d_i + 4\pi\sum_{i=2}^N\sum_{j=i+1}^N G_s(x_j;x_i) \right) + \pi\eps^2 \sum_{i=2}^N G_s(x_1;x_i) \right) + O(\eps^3)\,,
\eeq
and then upon rearranging terms we obtain
\beq
\ubar = \frac{\dfrac{\pi\eps}{4(N-1)} - \dfrac{\pi\eps^2}{N-1} \sum_{i=2}^N G_s(x_1;x_i)}{1 - \dfrac{2\eps}{\pi(N-1)} \left( \sum_{i=2}^N d_i + 4\pi\sum_{i=2}^N\sum_{j=i+1}^N G_s(x_j;x_i) \right)}\,,
\eeq
which can be further expanded into
\beq\label{eq:ubarapprox}
\ubar = \left( \frac{\pi\eps}{4(N-1)} - \frac{\pi\eps^2}{N-1} \sum_{i=2}^N G_s(x_1;x_i) \right) \left( 1 - \frac{2\eps}{\pi(N-1)} \left( \sum_{i=2}^N d_i + 4\pi\sum_{i=2}^N\sum_{j=i+1}^N G_s(x_j;x_i) \right) + O(\eps^2) \right)\,.
\eeq
Finally, by dropping all the $O(\eps^3)$ terms in \eqref{eq:ubarapprox} we find that $\ubar$ satisfies
\beq
\ubar = \frac{\pi\eps}{4(N-1)} + \frac{\eps^2}{2(N-1)^2}\sum_{i=2}^N d_i + \frac{2\pi\eps^2}{(N-1)^2}\sum_{i=2}^N\sum_{j=i+1}^N G_s(x_j;x_i) - \frac{\pi\eps^2}{N-1}\sum_{i=2}^N G_s(x_1;x_i) + O(\eps^3) \,.
\eeq

Next, we calculate the constants $C_j$ with $j=2,\,\ldots,\,N$ that control the outfluxes at each window. From eq.~\eqref{eq:twoTerms} we obtain that
\beq\label{eq:cjtemp}
C_j = \ubar\left( -\frac{2}{\pi} + \frac{4\eps}{\pi^2} \left(d_j + \sum_{\substack{i=2\\ i\neq j}}^N G_s(x_i;x_j)\right) \right) - 2\eps^2 G_s(x_1;x_j) + O(\eps^3)\,,
\eeq
and then by substituting the expression for $\ubar$ within \eqref{eq:cjtemp} and dropping higher order cubic terms yields
\begin{align*}
C_j &= -\frac{\eps}{2(N-1)} + \frac{2\eps^2}{N-1} \sum_{i=2}^N G_s(x_1;x_i) - 2\eps^2 G_s(x_1;x_j) - \frac{\eps^2}{\pi(N-1)^2} \sum_{i=2}^N d_i + \frac{\eps^2}{\pi(N-1)}d_j \\
&- \frac{4\eps^2}{(N-1)^2}\sum_{i=2}^N\sum_{k=i+1}^N G_s(x_k;x_i) + \frac{2\eps^2}{N-1}\sum_{\substack{i=2\\ i\neq j}}^N G_s(x_i;x_j) + O(\eps^3)\,.
\end{align*}
Next, by grouping all quadratic terms together and simplifying the summation terms we get
\begin{align}
C_j &= -\frac{\eps}{2(N-1)} - \frac{\eps^2}{\pi(N-1)^2}\left( \sum_{\substack{i=2\\ i\neq j}}^N \left(d_i - d_j\right) + 2\pi(N-1)\sum_{\substack{i=2\\ i\neq j}}^N \left( G_s(x_1;x_j) - G_s(x_1;x_i) - G_s(x_i;x_j)\right) \right. \nonumber \\
& \left. + 4\pi\sum_{i=2}^N\sum_{k=i+1}^N G_s(x_k;x_i) \right) + O(\eps^3)\,. \label{eq:cj}
\end{align}
At leading-order, we found that the outflux is the same for all exiting windows. The organization of the windows holes influences only via the higher order terms which depend on the Neumann-Green's function. However, for an arbitrary closed and bounded domain $\Omega$, the logarithmic singularity provides a correction that depends on the mean curvature function $H(x)$. Indeed, using the expression
\beq
d_i = -\frac{H(x_i)}{2} \log(\eps) + O(1)\,, \quad \text{for} \quad i \neq 1\,,
\eeq
we  find that \eqref{eq:cj} becomes
\beq\label{eq:cjh}
C_j = -\frac{\eps}{2(N-1)} + \frac{\eps^2\log(\eps)}{2\pi (N-1)^2} \sum_{\substack{i=2\\ i\neq j}}^N \left(H(x_i) - H(x_j)\right) + O(\eps^2)\,.
\eeq
This correction term vanishes for a boundary $\p\Omega$ with constant mean curvature (see section \S \ref{sec:sphere} for the unit sphere). We can now use our expressions for the average concentration $\ubar$ and for the constants $C_j$ controlling the outfluxes at each window to calculate the steady-state concentration value at the center of the patch $\p \Omega_{\eps_1}$ receiving an influx current. Upon setting $y = x_1$ in \eqref{eq:ubdyN}, we find
\beq
u(x_1) = \ubar + \int_{\p \Omega_{\eps_1}}G_s(x;x_1) dx + \sum_{j=2}^N \int_{\p \Omega_{\eps_j}}G_s(x;x_1)\frac{\p u}{\p n} dx \,,
\eeq
and then by using the expansion given in \eqref{eq:intref} for the integral of the Green's function near a singularity, we get
\beq
u(x_1) = \ubar + \eps - \frac{H(x_1)}{4} \eps^2\log(\eps) + \frac{H(x_1)}{8}\eps^2 + v_1 \pi \eps^2 + 2\pi\eps\sum_{j=2}^N C_jG_s(x_j;x_1)\,.
\eeq
Because we seek a general three-term asymptotic approximation for the concentration drop, we only keep the leading-order term of the expression \eqref{eq:cjh} for the constants $C_j$, and this yields
\beq
u(x_1) = \ubar + \eps - \frac{H(x_1)}{4} \eps^2\log(\eps) + \frac{H(x_1)}{8}\eps^2 + v_1 \pi \eps^2 - \frac{\pi\eps^2}{N-1}\sum_{j=2}^N G_s(x_j;x_1) + O(\eps^3)\,,
\eeq
and by substituting the expansion for $\ubar$ we get,
\beq
\begin{split}\label{eq:ux1temp}
u(x_1) &= \left( 1 + \frac{\pi}{4(N-1)}\right)\eps - \frac{H(x_1)}{4} \eps^2\log(\eps) + \frac{H(x_1)}{8}\eps^2 + v_1 \pi \eps^2 - \frac{\pi\eps^2}{N-1}\sum_{j=2}^N G_s(x_j;x_1) \\
& + \frac{\eps^2}{2(N-1)^2}\sum_{i=2}^N d_i + \frac{2\pi\eps^2}{(N-1)^2}\sum_{i=2}^N\sum_{j=i+1}^N G_s(x_i;x_j) - \frac{\pi\eps^2}{N-1}\sum_{i=2}^N G_s(x_1;x_i) + O(\eps^3)\,.
\end{split}
\eeq
Using the expression for $d_i$ found in \eqref{eq:d_i} we compute that
\beq
\frac{\eps^2}{2(N-1)^2}\sum_{i=2}^N d_i = -\frac{\eps^2\log(\eps)}{2(N-1)^2}\sum_{i=2}^N H(x_i) + \frac{1-\log(2)\eps^2}{4(N-1)^2}\sum_{i=2}^N H(x_i) + \frac{\pi\eps^2}{(N-1)^2}\sum_{i=2}^N v_i\,,
\eeq
and by collecting the various terms in $\eps$, $\eps^2\log(\eps)$ and $\eps^2$ in \eqref{eq:ux1temp}, we obtain the asymptotic expansion given below,
\beq\label{eq:ux1}
\begin{split}
u(x_1) &= \left( 1 + \frac{\pi}{4(N-1)}\right)\eps - \left(H(x_1) +  \frac{1}{(N-1)^2}\sum_{i=2}^N H(x_i)\right)\frac{\eps^2}{4}\log(\eps) \\
&+ \left(\frac{H(x_1)}{8} + \frac{1-\log(2)}{4(N-1)^2}\sum_{i=2}^N H(x_i) + \pi v_1 + \frac{\pi}{(N-1)^2} \sum_{i=2}^N v_i \right. \\
& \left. + \frac{2\pi}{(N-1)^2}\sum_{i=2}^N\sum_{j=i+1}^N G_s(x_i;x_j) - \frac{2\pi}{N-1}\sum_{i=2}^N G_s(x_1;x_i) \right)\eps^2 + O(\eps^3)\,.
\end{split}
\eeq
As in \S \ref{sec:nbc}, the first two terms of the asymptotic approximation \eqref{eq:ux1} depend on the mean curvature at the center of the patches $\p\Omega_{\eps_j}$ and are independent of the distances $\|x_i - x_j\|$. Hence, by dropping the third quadratic term, we obtain a two-term asymptotic expansion for the concentration drop resulting from influx diffusion, valid for arbitrary domains with smooth boundaries. In summary, we have
\begin{res}
For $N$ well-separated circular windows $\p\Omega_{\eps_j}$ of small radius $\eps \ll 1$, all absorbing except for $j=1$, which receives an influx current, the normalized solution $u(x)$ of \eqref{eq:pdenonN} in the domain $\Omega$ has a two-term asymptotic approximation for the concentration drop given by
\beq
u(x_1) = \eps\left(1+\frac{\pi}{4(N-1)} - \left(H(x_1) + \frac{1}{(N-1)^2}\sum_{i=2}^N H(x_i)\right)\frac{\eps}{4}\log(\eps) + O(\eps)\right)\,.
\eeq
where $H(x_j)$ for $j=1,\ldots,N$ are the mean curvature at $x_j$ on the boundary $\p\Omega$. For the special case of $N=2$ patches, the formula reduces to
\beq
u(x_1) = \eps\left(1+\frac{\pi}{4} - \frac{H(x_1) + H(x_2)}{4}\eps\log(\eps) + O(\eps)\right)\,.
\eeq
Furthermore, the outward flux $\Phi_j$ at each exiting window $\p \Omega_{\eps_j}$ $j=2,\ldots,N$ is approximated by
\beq
\Phi_j = 2\pi C_j\eps = -\frac{\pi\eps^2}{(N-1)} + \frac{\eps^3\log(\eps)}{(N-1)^2} \sum_{\substack{i=2\\ i\neq j}}^N \left(H(x_i) - H(x_j)\right) + O(\eps^3) \quad j \neq 1\,.
\eeq
\end{res}
\section{Unit Sphere Domain}\label{sec:sphere}
In this section, we provide explicit three-term asymptotic formulas for the normalized concentration difference on a unit sphere domain $\Omega = \{x|\|x\| \leq 1\}$, that we compare against numerical solutions \cite{comsol}. Earlier, we found an exact formula for a general surface using the Neumann-Green's function. 
\subsection{Asymptotic expression for the concentration in a unit ball}
For the unit sphere, the Green's function $G_s(x;y)$ solution of \eqref{eq:greenfun} \cite{cheviakov2010}  is given explicitly (with a zero mean) by
\beq\label{eq:gsphere1}
G_s(x;y) = \frac{1}{2\pi\|x-y\|} + \frac{1}{8\pi}\left(\|x\|^2 + 1\right) + \frac{1}{4\pi}\log\left(\frac{2}{1-\|x\|\cos(\gamma)+\|x-y\|}\right) - \frac{7}{10\pi}\,, \quad \|y\| = 1\,,
\eeq
where $\gamma$ is the angle between $x$ and $y$. In order to evaluate formulas \eqref{eq:NeumannSolDiff} and \eqref{eq:ux1}, we shall compute $G_s(x;y)$ for $x$ on the boundary $\p\Omega$. By setting $\|x\| = 1$ in \eqref{eq:gsphere1} and using that $\cos(\gamma) = 1 - \|x-y\|^2/2$ from the cosine law, we find that
\beq\label{eq:gsphere2}
G_s(x;y) = \frac{1}{2\pi\|x-y\|} - \frac{1}{4\pi}\log\left(\frac{1}{2}\|x-y\|^2+\|x-y\|\right) + \frac{\log(2)}{4\pi} - \frac{9}{20\pi}\,, \quad \|x\| = \|y\| = 1,\, x \neq y \,.
\eeq
To recover the expansion given in \eqref{eq:gdec} we use that $r = \|x-y\|$ in \eqref{eq:gsphere2} and we further obtain
\beq\label{eq:gsphere3}
G_s(x;y)  = \frac{1}{2\pi r} - \frac{1}{4\pi}\log\left(\frac{r^2}{2} + r\right) + \frac{\log(2)}{4\pi} - \frac{9}{20\pi}.
\eeq
Finally, we expand the logarithmic singularity for $r$ small 
\beq\label{eq:gsphere4}
G_s(x;y) \sim \underbrace{\frac{1}{2\pi r} - \frac{1}{4\pi}\log\left(r\right)}_{g(r)} + \underbrace{\frac{\log(2)}{4\pi} - \frac{9}{20\pi}}_{v} + O(r) \,,
\eeq
where $g(r)$ is the sum of the coulomb and logarithmic singularities while $v$ stands for the regular part of the Green's function. For the unit sphere, the mean curvature on the boundary is constant, $H \equiv 1$ and near the singularity, the integral of the $O(r)$ term using polar coordinates behaves like
\beq
2\pi\int_0^\eps O(r)rdr \sim O(\eps^3)\,,
\eeq
and can thus be neglected. The present asymptotic formulas for the steady-state concentration drop due to influx diffusion on the unit sphere and for an exiting flux at each absorbing window are
summarized as follows:
\begin{pres}\label{pres:pres1}
For the unit ball $\Omega = \{x|\,\|x\| \leq 1\}$ with two well-separated circular patches of radius $\eps \ll 1$ on its spherical boundary, with one of them receiving while the other is emitting a normalized influx current, a three-term asymptotic solution for the normalized concentration drop is given by
\beq\label{eq:sphere2}
u(x_1) - u(x_2) = 2\eps - \frac{\eps^2}{2}\log\left(\eps\right) + \left(\frac{1}{4} - \frac{1}{l} + \frac{1}{2}\log\left(\frac{l^2}{2} + l\right) \right)\eps^2 + O(\eps^3)\,,
\eeq
where $l = \|x_1 - x_2\| \sim O(1)$ is the distance between the centers of each patch. Alternatively if there are $N$ identical, well-separated, circular patches on the boundary, all absorbing with the exception of the first one which receives a normalized influx current, then the three-term asymptotic approximation for the normalized concentration drop is given by
\beq\label{eq:sphereN}
\begin{split}
u(x_1) &= \left(1 + \frac{\pi}{4(N-1)}\right)\eps - \frac{N}{4(N-1)}\eps^2\log(\eps) + \left(\frac{N+1}{8(N-1)} - \frac{\log(2)}{4(N-1)} \right. \\
& \left. + \frac{1}{(N-1)^2}\sum_{i=2}^N\sum_{j=i+1}^N \left( \frac{1}{l_{ij}} - \log\left(\frac{l_{ij}^2}{2} + l_{ij}\right) \right) - \frac{1}{N-1}\sum_{i=2}^N \left( \frac{1}{l_{i1}} - \frac{1}{2}\log\left(\frac{l_{i1}^2}{2} + l_{i1}\right) \right)\right)\eps^2 + O(\eps^3)\,,
\end{split}
\eeq
where $l_{ij} = \|x_i - x_j\| \sim O(1)$ corresponds to the different distances. This formula reduces to
\beq\label{eq:sphere2abs}
u(x_1) = \left(1+\frac{\pi}{4}\right) \eps - \frac{\eps^2}{2}\log(\eps) + \left( \frac{3}{8} - \frac{\log(2)}{4} - \frac{1}{l} + \frac{1}{2}\log\left( \frac{l^2}{2} + l \right) \right)\eps^2 + O(\eps^3)\,,
\eeq
when there are $N=2$ circular patches on the boundary $\p \Omega$.

Finally, we recall that the total flux at each exiting window $\p\Omega_{\eps_j}$ $j=2,\ldots,N$ is given by
\beq
\phi_j = \int_{\p\Omega_{\eps_j}} \frac{\p u(x)}{\p n} dx = 2\pi C_j\eps\,, \quad j \neq 1\,,
\eeq
and thus using the expansion for the constants $C_j$ given in \eqref{eq:cj} and from the fact that the mean curvature and the regular part of the Green's function are constant on the sphere, we have
\beq
\begin{split}
C_j &= -\frac{\eps}{2(N-1)} - \left(\sum_{\substack{i=2\\ i\neq j}}^N\left(g(l_{1j}) - g(l_{1i}) - g(l_{ij}) \right) \right. \left. + \frac{2}{(N-1)} \sum_{i=2}^N \sum_{k=i+1}^N g(l_{ik})\right)\frac{2\eps^2}{(N-1)} + O(\eps^3)\,.
\end{split}
\eeq
For the total flux $\Phi_j$, we therefore obtain the two-term expansion
\beq\label{eq:fluxSphere}
\begin{split}
\phi_j &= -\frac{\pi\eps^2}{N-1} - \left(\sum_{\substack{i=2\\ i\neq j}}^N\left(\frac{1}{l_{1j}} - \frac{1}{l_{1i}} - \frac{1}{l_{ij}} - \frac{1}{2}\log\left( \frac{2l_{1j}^2 + 4l_{1j}}{\left(l_{1i}^2 + 2l_{1i}\right)\left(l_{ij}^2 + 2l_{ij}\right)} \right) \right) \right. \\
& \left. + \frac{2}{(N-1)} \sum_{i=2}^N \sum_{k=i+1}^N \left(\frac{1}{l_{ik}} - \frac{1}{2}\log\left( \frac{l_{ik}^2}{2} + l_{ik} \right) \right)\right)\frac{2\eps^3}{(N-1)} + O(\eps^4)\,, \quad j \neq 1\,.
\end{split}
\eeq
\end{pres}
We find that the concentration difference is independent of the regular part of the Green's function. This is obvious when there are only $N=2$ circular patches, but for general $N$ it follows from the fact that
\begin{equation*}
\pi v\left(1 + \frac{1}{(N-1)^2}\sum_{i=2}^N 1 + \frac{2}{(N-1)^2} \sum_{i=2}^N\sum_{j=i+1}^N 1 - \frac{2}{N-1}\sum_{i=2}^N 1 \right) = 0\,.
\end{equation*}
\subsection{Comparing asymptotic vs numerical simulations}
In this section, we compare various asymptotic formula to numerical simulations for various window configurations on the unit ball. We study how the radius of each window and the distances between them affect the concentration drop. We will use the asymptotic formulas \eqref{eq:sphere2},\eqref{eq:sphere2abs} and \eqref{eq:sphereN} against numerical simulations performed with COMSOL \cite{comsol}, where the center of the window $\p\Omega_{\eps_1}$ receiving an influx current is conveniently located at the North pole. To measure the discrepancy between the asymptotic and numerical solutions, we compute the relative error
\beq
Re = 100 \times \frac{\delta u_{\rm asym} - \delta u_{\rm num}}{\delta u_{\rm num}}\,,
\eeq
where $\delta u \equiv u(x_1) - u(x_j)$ with $j\neq 1$. Our results are presented in Fig.~\ref{fig:fig2} - \ref{fig:fixed_area} below. The best agreement between the asymptotic and the numerics is obtained for well-separated as opposed to closely located windows.\\
Qualitatively similar results are obtained in Fig.~\ref{fig:fig2} and \ref{fig:fig3} for $N=2$ holes. For a fixed radius, we observe a weak interaction between the two windows when they are far apart. Furthermore, smaller concentration differences are obtained when the exiting window has an absorbing boundary condition, as opposed to prescribing the outward flux via a Neumann boundary condition. This follows from the leading-order term of \eqref{eq:sphere2} being bigger than in \eqref{eq:sphere2abs}, as seen from the inequality $\left(1+\frac{\pi}{4}\right)\eps < 2\eps$.
\begin{figure}[H]
\includegraphics[width=\linewidth]{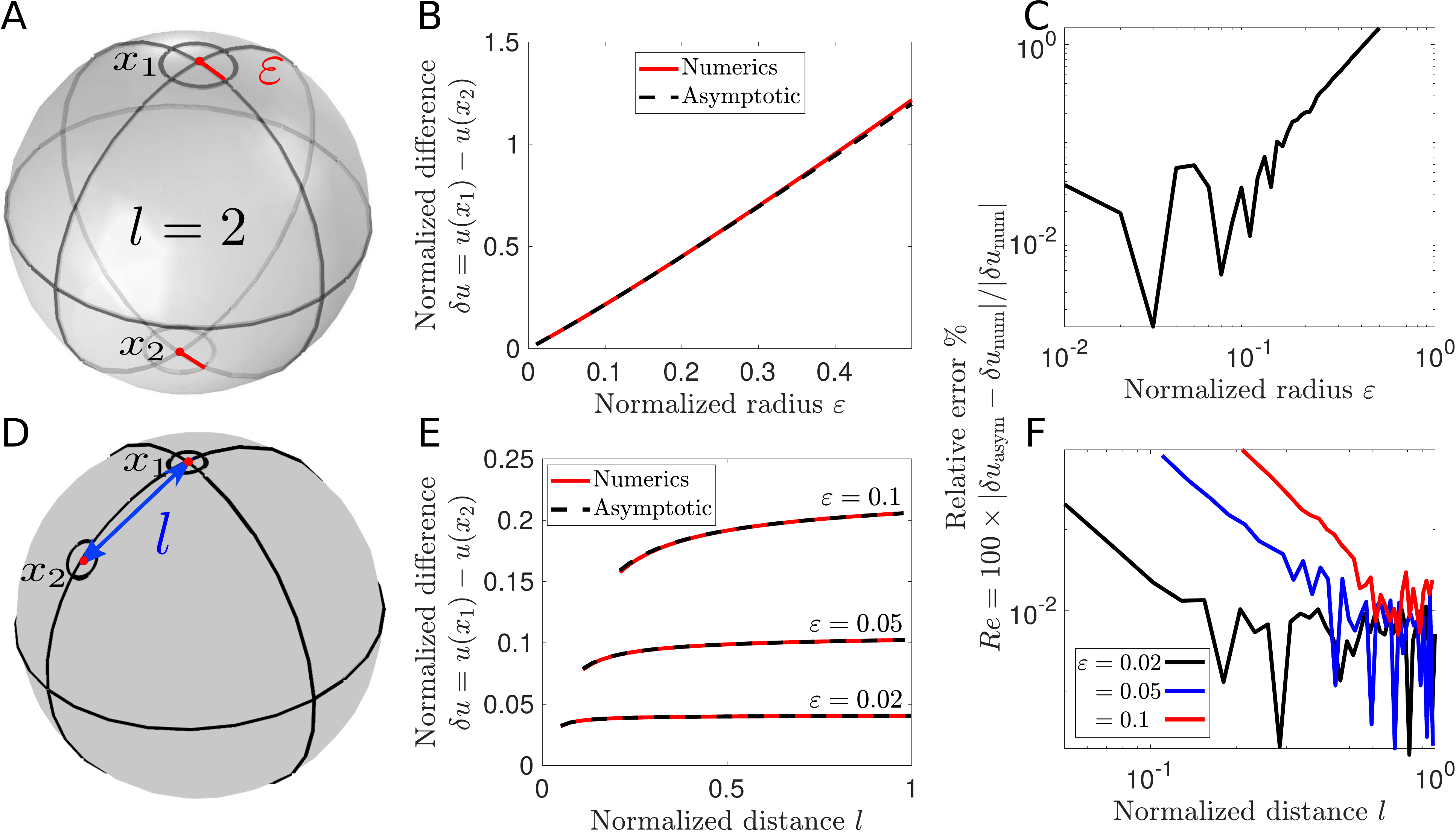}
\caption{\textbf{Concentration drop in a unit ball with two windows with Neumann boundary conditions (formula \eqref{eq:sphere2}).} \textbf{A}: An influx and an outflux currents are applied respectively on two patches centered around the North $(x_1)$ and South $(x_2)$ Poles. \textbf{B}: Concentration difference as a function of the radius $\eps$. \textbf{C}: Relative error as a function of the radius $\eps$. \textbf{D}: By increasing the colatitude of a patch centered in $x_2$ on the trivial azymuth, we study the effect of varying the distance $l = \|x_1-x_2\|$. \textbf{E}: Concentration difference as a function of the distance $l$ for $\eps = 0.02,\,0.05,\, 0.1$. \textbf{F}: Relative error as a function of the distance $l$.}
\label{fig:fig2}
\end{figure}
\begin{figure}[H]
\includegraphics[width=\linewidth]{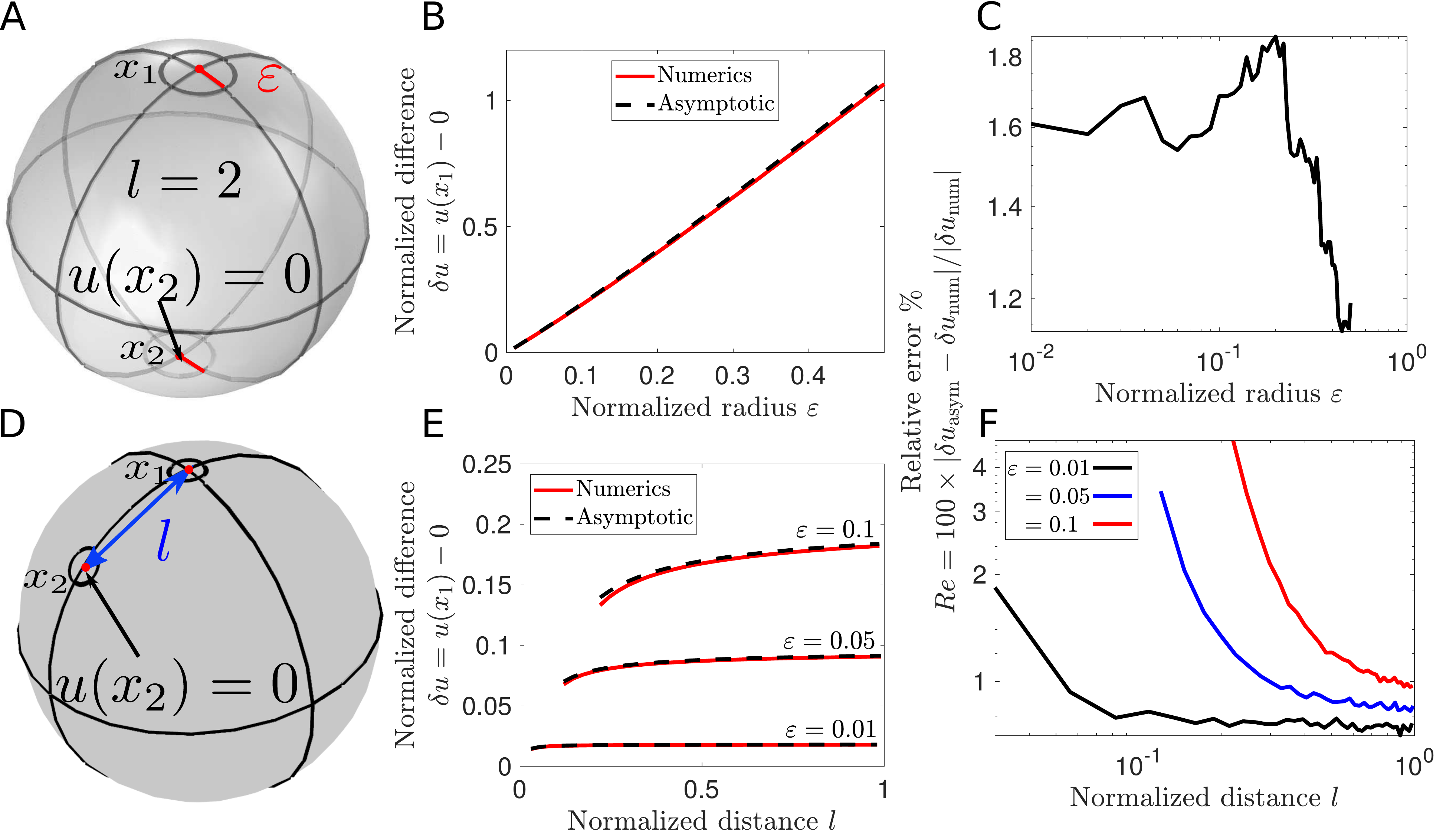}
\caption{\textbf{Concentration drop in a unit ball with a single absorbing window (formula \eqref{eq:sphere2abs}).} \textbf{A}: The absorbing window is centered at the South Pole $(x_2)$. \textbf{B}: Concentration difference as a function of the radius $\eps$. \textbf{C}: Relative error as a function of the radius $\eps$. \textbf{D}: By increasing the colatitude of a patch centered in $x_2$ on the trivial azymuth, we study the effect of varying the distance $l = \|x_1-x_2\|$. \textbf{E}: Concentration difference as a function of the distance $l$ for $\eps = 0.01,\,0.05,\, 0.1$. \textbf{F}: Relative error as a function of the distance $l$.}
\label{fig:fig3}
\end{figure}
For a configuration with several exiting holes, we found a good agreement between the asymptotic and numerical solutions (Fig.~\ref{fig:fig4}A-C), for either two, four or six windows located at the North and South poles, and equidistantly on the equator. A smaller concentration difference is to be expected as the number of holes increases. The pattern of holes is different in Fig.~\ref{fig:fig4}D-F, where we vary the radius of a concentric ring where the exiting windows are located. As the radius of the ring increase, the concentration drop slowly increases. We reported a similar relation with the distance in Fig.~\ref{fig:fig2} and \ref{fig:fig3} for $N=2$ windows.
\begin{figure}[H]
\centering
\includegraphics[width=\linewidth]{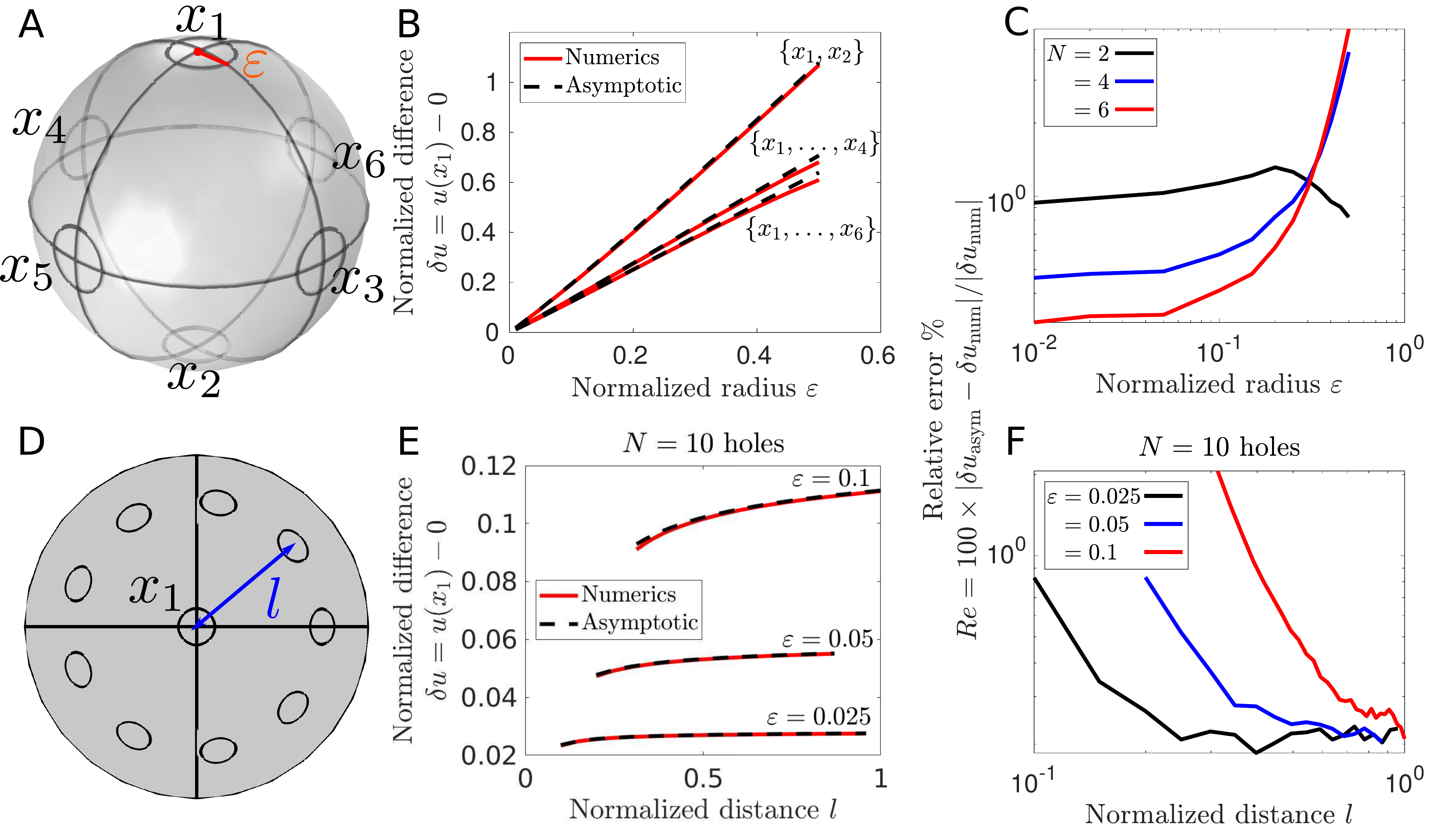}
\caption{\textbf{Several absorbing holes (formula \eqref{eq:sphereN}).} \textbf{A}: There are 6 equidistant holes on the sphere, located at the North and South Poles with the last four on the equator. \textbf{B}: Normalized concentration drop as a function of the radius $\eps$ for $N=2,\,4,\,6$ holes. \textbf{C}: Relative error as a function of $\eps$. \textbf{D}: We consider $9$ holes located on a concentric ring of radius $l$ centered at the North Pole (where the usual influx patch is located). \textbf{E}: Concentration drop as a function of the distance $l$. \textbf{F}: Relative error as a function of the distance $l$.}
\label{fig:fig4}
\end{figure}
We also explored different types of clustering on the concentration drop across the domain: for a lineic (Fig.~\ref{fig:figlineic}) window configuration and for exiting windows that form circular clusters (in Fig.~\ref{fig:figcluster}). We obtain similar results for both configurations and as previously observed, higher concentration differences are obtained when exiting windows move away from the influx location.  When exiting windows form clusters, the concentration drop across the domain tends to be slightly larger as shown in Fig.~\ref{fig:figlineic}D-F and \ref{fig:figcluster}D-F. This result suggests that absorbing windows are more effective in catching diffusing particles when they are well separated and widely distributed on the boundary.
\begin{figure}[H]
\centering
\includegraphics[width=\linewidth]{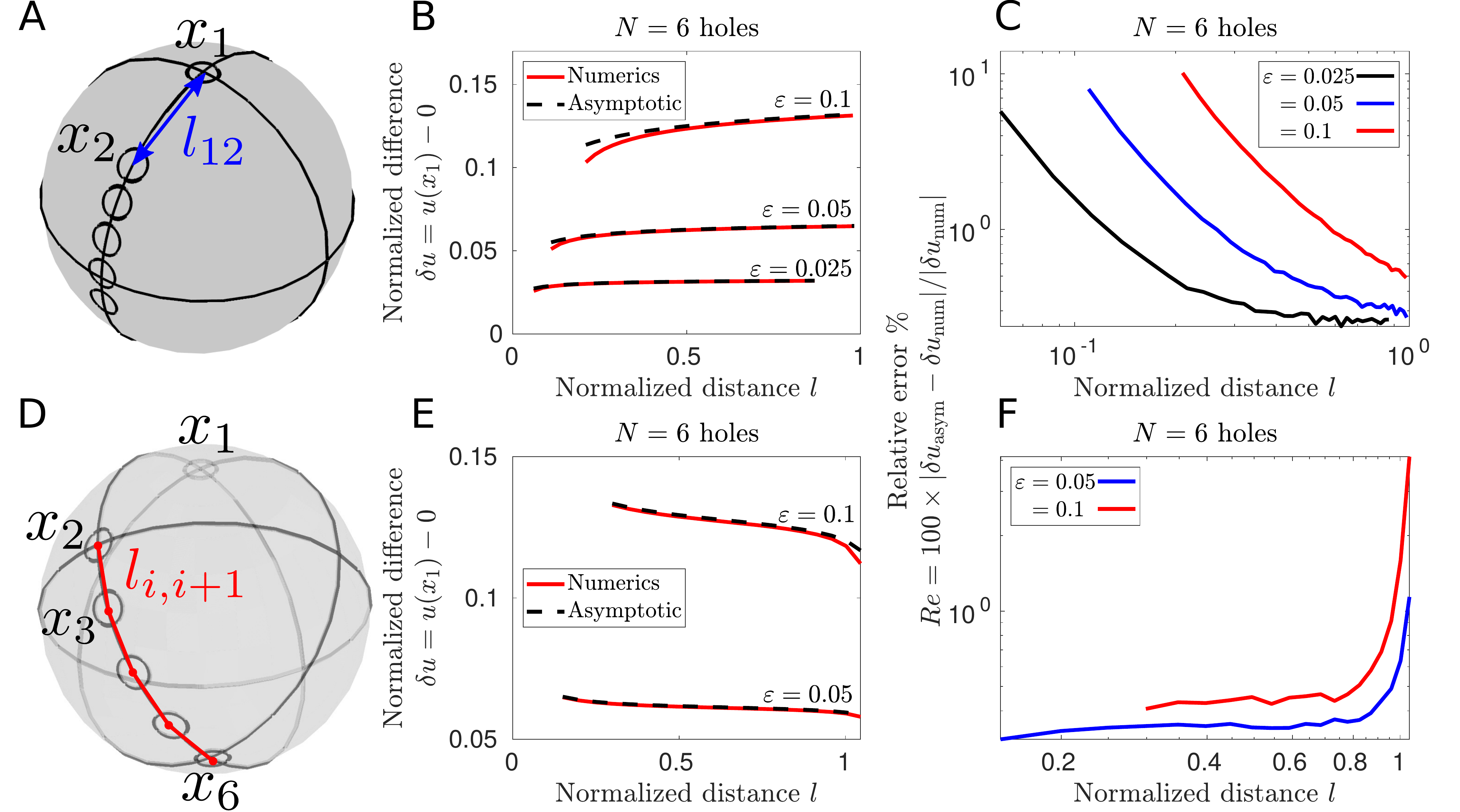}
\caption{\textbf{Concentration drop due to a linear configuration of absorbing windows (formula \eqref{eq:sphereN}).} \textbf{A}: There are $5$ absorbing windows located on the trivial azymuth, with a distance of $2.5\eps$ between each neighbor. \textbf{B}: Concentration drop as a function of the distance from the North pole to the center of the nearest window on the line. \textbf{C}: Relative error as a function of the distance $l$. \textbf{D}: Here $l$ represents the distance between each neighboring window. The center of the second window remains fixed on the Equator. \textbf{E}: Concentration drop as a function of $l$. \textbf{F}: Relative error as a function of $l$.}
\label{fig:figlineic}
\end{figure}
\begin{figure}[H]
\centering
\includegraphics[width=\linewidth]{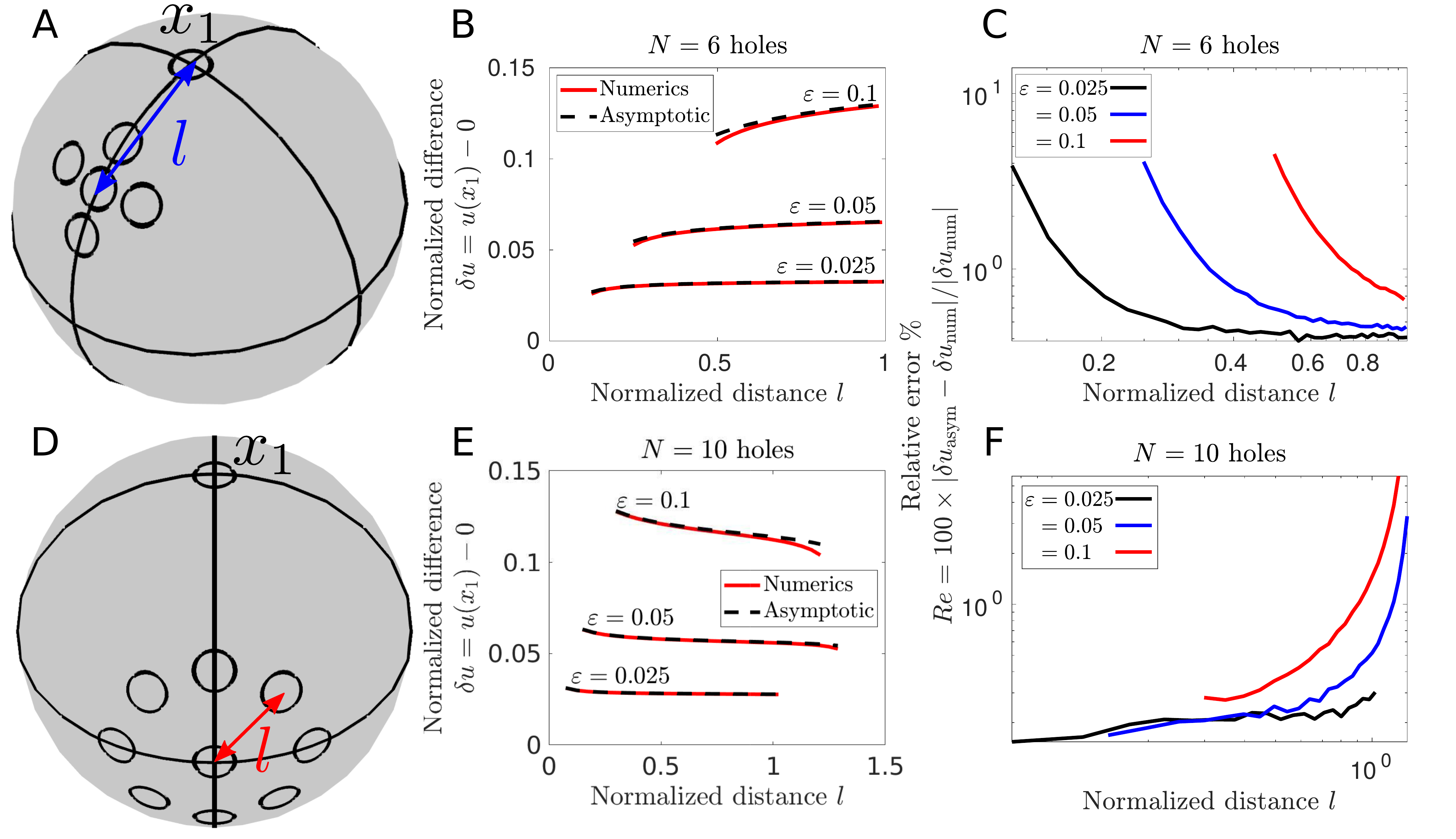}
\caption{\textbf{Concentration drop due to a cluster of absorbing windows (formula \eqref{eq:sphereN}).} \textbf{A}: Here $l$ represents the distance from the North Pole to the center of a cluster of radius $2.5\eps$. \textbf{B}: Concentration drop as a function of $l$. \textbf{C}: Relative error as a function of $l$. \textbf{D}: Here $l$ corresponds to the radius of a cluster centered on the Equator. \textbf{E}: Concentration drop as a function of $l$. \textbf{F}: Relative error as a function of $l$.}
\label{fig:figcluster}
\end{figure}
In Fig.~\ref{fig:fixed_area}, we increase the number of absorbing windows while keeping fixed the total window area to be $\frac{\pi}{5}$, including the influx receiving window, with all exiting windows located on the equator at equal distances from each other. To compensate the higher number of absorbing windows, the radius $\eps$ must decrease and thus lower concentration differences across the domain are observed.
\begin{figure}[H]
\centering
\includegraphics[width=\linewidth]{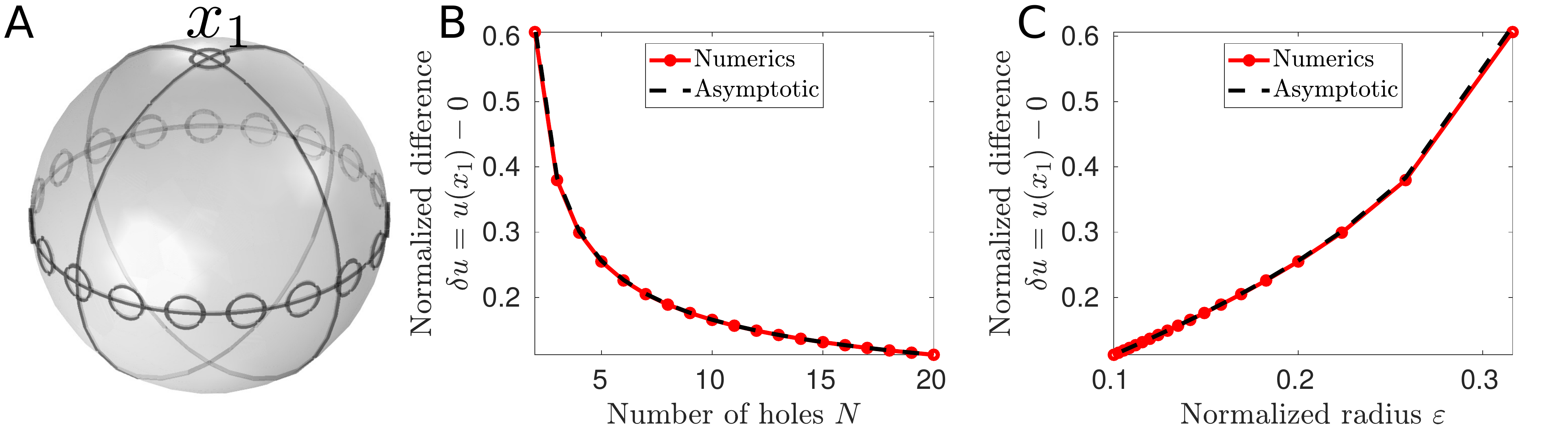}
\caption{\textbf{Variation of the number of absorbing windows (formula \eqref{eq:sphereN}).} \textbf{A}: The usual influx hole is located at the North Pole and the absorbing windows are equidistantly located along the Equator. \textbf{B}: Concentration drop as a function of the number of holes $N$, from $2$ to $20$. The total window area remains constant and equal to $\pi/5$. \textbf{C}: To compensate the addition of windows along the Equator their radius must decrease.}
\label{fig:fixed_area}
\end{figure}
Finally in Fig.~\ref{fig:weberFlux}, we show that the exiting fluxes are affected by the geometrical organization of the absorbing windows: in that case, the total flux magnitude is higher for holes that are located near the window receiving the influx.
\begin{figure}[H]
\centering
\includegraphics[width=\linewidth]{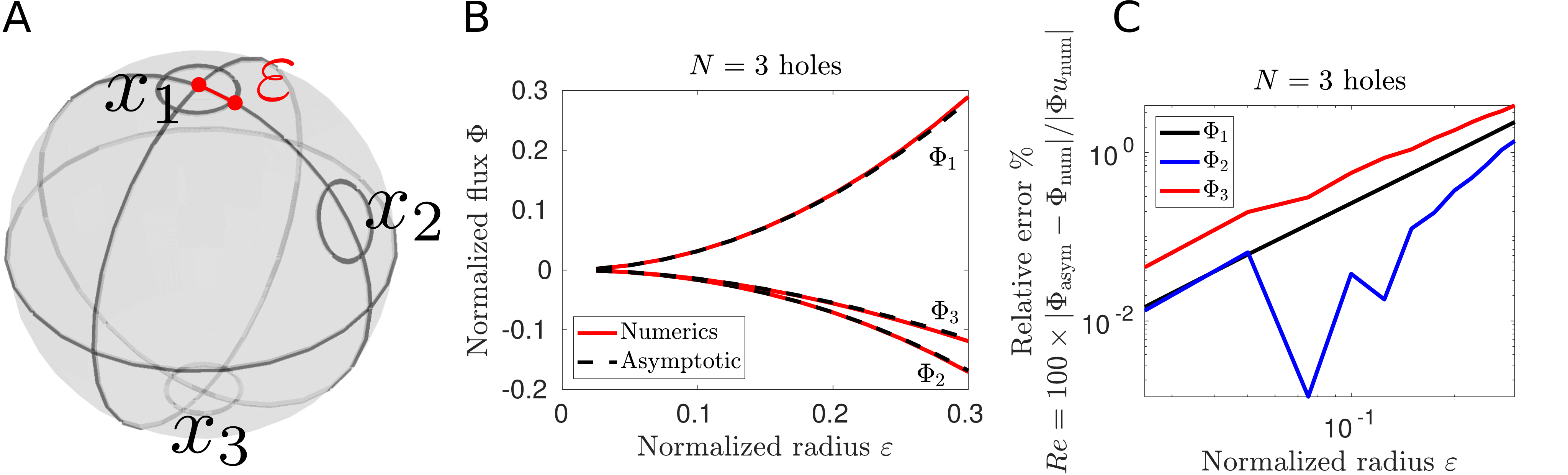}
\caption{\textbf{Total influx versus exiting fluxes (formula \eqref{eq:fluxSphere}).} \textbf{A}: The influx window is located at the North Pole, the second window has colatitude $\frac{\pi}{3}$, while the third one is centered at the South Pole. \textbf{B}: Total window fluxes vs the radius $\eps$. \textbf{C}: Relative error between the asymptotic formula \eqref{eq:fluxSphere} and numerically computed fluxes for each window.}
\label{fig:weberFlux}
\end{figure}
\section{Penetration length: a spatial scale to quantify the flow between two neighboring windows using}\label{sec:plength}
\begin{figure}[H]
\centering
\includegraphics[width=\linewidth]{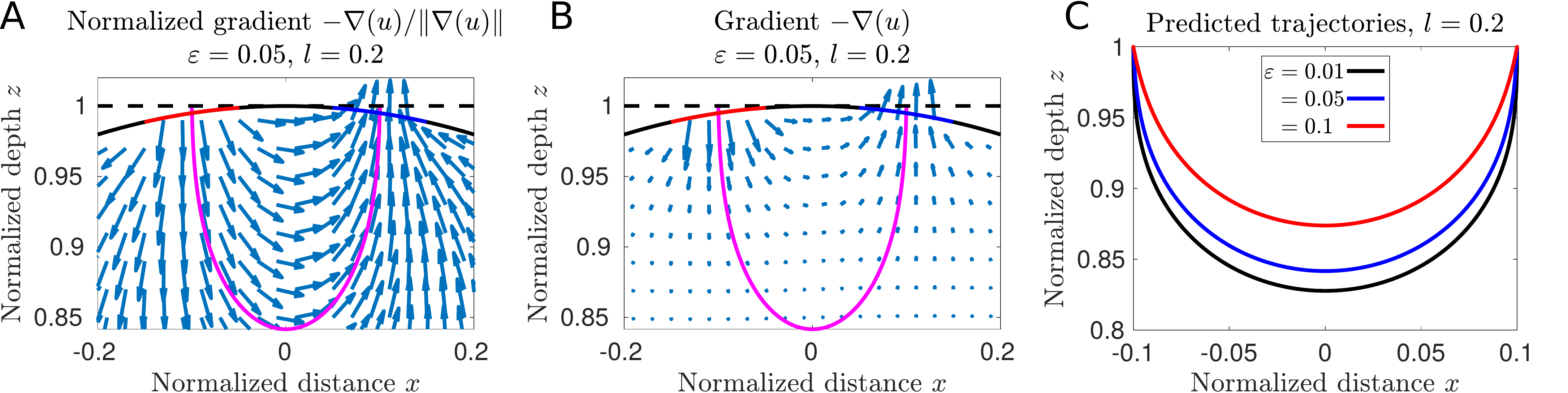}
\caption{\textbf{Flow lines and influx trajectories}. \textbf{A}: A trajectory (purple) from a boundary follows the concentration gradient vector field (normalized blue arrow) for a distance $l=0.2$ between the windows of size $\eps=0.05$. \textbf{B}: Same as Panel A for a non-normalized gradient vector field. The arrow length is proportional to the magnitude of the gradient. \textbf{C}: Examples of 3 different trajectories associated to different radius $\eps=0.01$ (black) $0.05$ (blue) and $0.1$ (red), for a planar boundary with a distance $l=0.2$ between the start and end points. As the window radius $\eps$ increases the length of penetration decays.}
\label{fig:fig9}
\end{figure}

In this final section, to estimate how deep inside a medium a concentration change can propagate following an influx from one window and an outflux from a second one, we introduce a length scale that we call the penetration length $L_{pe}$. We evaluate how this length depends on the distance between the windows, the current $I$ of particles entering through window 1 and the diffusion coefficient $D$ of the moving particles. The penetration length is defined for a given window configuration and an influx $I$ as the maximal distance to the boundary of the trajectory associated with the steepest concentration descent starting at the center of window 1 and exiting at the center of window 2. When the two windows are located at the North and South poles respectively of a ball of radius $R$, then $L_{pe}=2R$ and the length is independent of the influx I and the radius $\eps$ of the two windows. However, in general it is not clear how this distance depends on the local geometry of the two windows. To estimate how $L_{pe}$ depends on various parameters, we employ an explicit solution of Laplace's equation. \\
We focus on the half-space domain $\Omega = \{(x,y,z) \, \in \R^3| \, z \geq 0\}$, with two circular patches of radius $\eps$ on its planar boundary $\p \Omega = \{(x,y,z) \, \in \R^3 | z = 0\}$  separated by a distance $l$, and centered at $\left(-\frac{l}{2},0,0\right)$ and $\left(\frac{l}{2},0,0\right)$ respectively. Then the steady-state solution $u(x,y,z)$ of the diffusion equation \eqref{eq:sspde} when an influx $I$ is fixed on window one is obtained by integrating Laplace's equation
\beq
\frac{\p ^2 u}{\p  x^2} + \frac{\p ^2 u}{\p  y^2} + \frac{\p ^2 u}{\p  z^2} = 0\,,
\eeq
with mixed boundary conditions
\beq
\left.-\frac{\p  u}{\p  z}\right|_{z=0} =
\begin{cases}
I & \sqrt{\left(x+\frac{l}{2}\right)^2 + y^2} \leq \eps \\ &\\
-I & \sqrt{\left(x-\frac{l}{2}\right)^2 + y^2} \leq \eps \\&\\
0 & \text{elsewhere}
\end{cases}\,.
\eeq
The solution has an integral representation (see appendix)
\beq\label{eq:uplane}
u(x,y,z) = u_0  + \eps I \int_0^\infty e^{-mz}J_1(m\eps)\left( J_0\left(m\sqrt{\left(x+\frac{l}{2}\right)^2 + y^2}\right) - J_0\left(m\sqrt{\left(x-\frac{l}{2}\right)^2 + y^2}\right) \right) \frac{dm}{m}\,,
\eeq
where $u_0$ is an arbitrary positive constant and $J_n(x)$ is the Bessel function of order $n$ \cite{Abramowitz}. We study now the trajectory of a flux of particles starting at the center of the influx window $\left(-\frac{l}{2},0,0\right)$. A trajectory is shown in Fig.~\ref{fig:fig9}A-B with the associated vector field. Three trajectories are represented in Fig.~\ref{fig:fig9}C, showing that the penetration length $L_{pe}$ decays with the radius $\eps$.
We start by generating a 3-D curve $(x(t),y(t),z(t))$ tangent to the concentration gradient $\nabla u$, and thus is solution of
\beq
\frac{d}{dt} \begin{pmatrix} x(t) \\ y(t) \\ z(t) \end{pmatrix} = - \nabla u\left(x(t),y(t),z(t)\right)\,,
\eeq
where the minus sign imposes that the gradient vector field points from high to low concentration areas. The center of the exiting window is $(\frac{l}{2},0,0)$, and by symmetry the trajectory lies in the $y=0$ plane. A parametrization of the 2-D curve $C(t) = (x(t),z(t))$ starting at $C(0) = \left(-\frac{l}{2},0\right)$ and tangent to the vector field satisfies
\beq \label{penel}
\frac{dC(t)}{dt} = - \left.
\begin{pmatrix}
{\frac{\p  u}{\p  x}} \\ \\
{\frac{\p  u}{\p  z}}
\end{pmatrix}\right|_{(x(t),0,z(t))}\,,
\eeq
which can be written using eq.~\eqref{eq:uplane} as an integro-differential equation
\beq\label{eq:ode}
\begin{pmatrix}
\dot{x}(t) \\ \\
\dot{z}(t)
\end{pmatrix}
= \eps I
\begin{pmatrix}
\int_0^\infty e^{-mz}J_1(m\eps)\left(J_1\left(m\left(\frac{l}{2} + x\right)\right) + J_1\left(m\left(\frac{l}{2} - x\right)\right)\right)dm \\ \\
\int_0^\infty e^{-mz}J_1(m\eps)\left(J_0\left(m\left(\frac{l}{2} + x\right)\right) - J_0\left(m\left(\frac{l}{2} - x\right)\right)\right)dm \\
\end{pmatrix}\,.
\eeq
The penetration length $L_{pe}$ and the travel time duration $T_{tr}$ between the two centers are defined by
\beq
L_{pe}= \max_{t > 0} z(t)\,, \quad T_{tr} = \inf \left\{t>0 \, \left| \, C(t) = \left(\frac{l}{2},0\right) \quad \text{given that} \quad C(0) = \left(-\frac{l}{2},0 \right) \right. \right\}\,.
\eeq
By definition, and using the symmetry of the domain, the penetration length $L_{pe}$ is achieved at equal distance of the two windows at $x=0$, for exactly half of the travel time, and thus we have $C(T_{tr}/2) = (0,L_{pe})$. At this point, the tangent vector is parallel to the $z = 0$ plane, and thus $\dot{z}(T_{tr}/2) = 0$. Furthermore, the amplitude of the current only affects the time-scale of the trajectory, and not the penetration length, as it can absorbed by a change of time in system \eqref{eq:ode}. However, for an influx with weak amplitude $I$, we expect longer travel times $T_{tr}$, and vice-versa for larger amplitudes $I$.

\begin{figure}[H]
\centering
\includegraphics[width=0.66\linewidth]{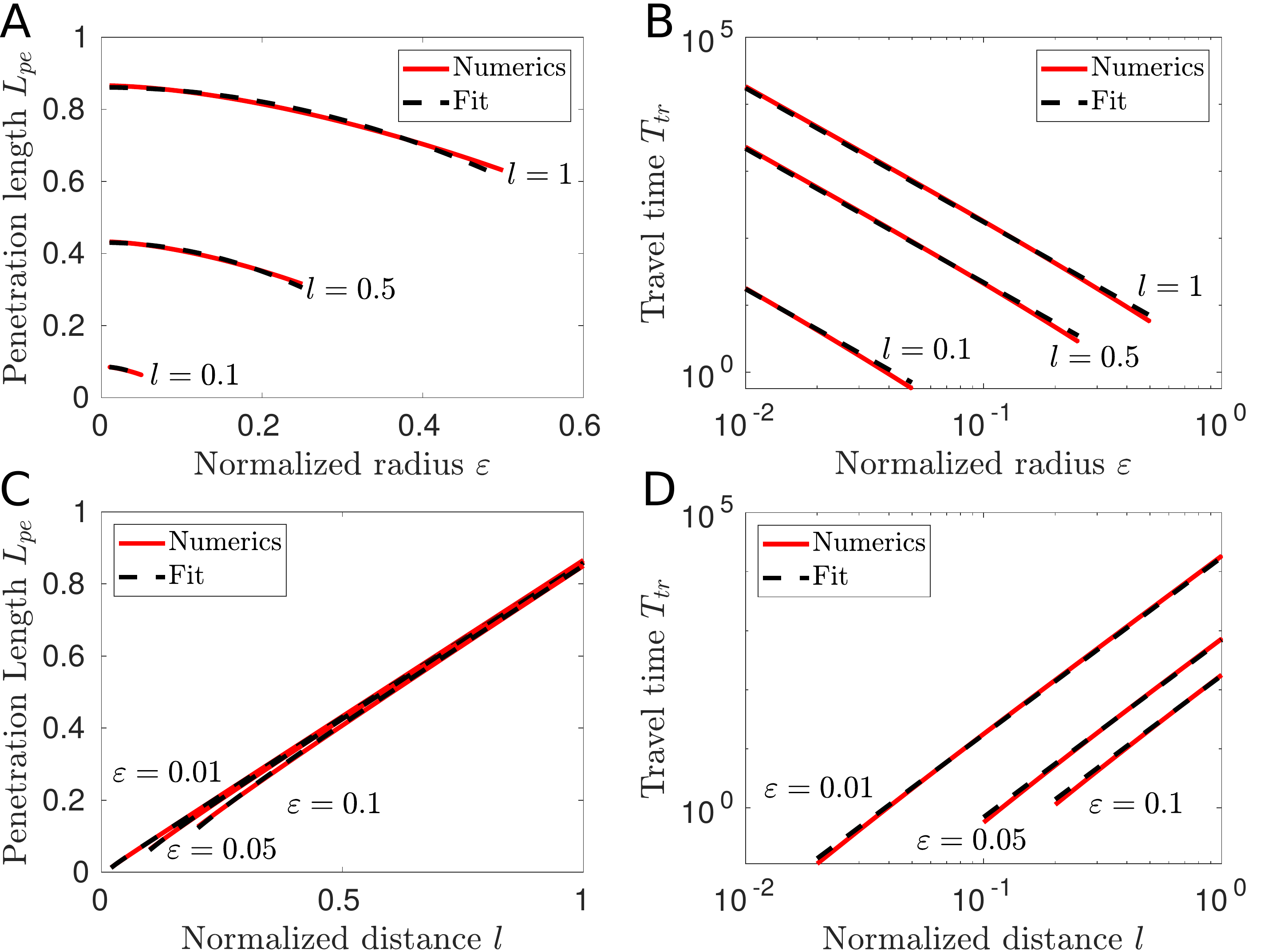}
\caption{\textbf{Properties of penetration length $L_{pe}$ and travel time $T_{tr}$.} \textbf{A}: $L_{pe}$ decays with the window radius $\eps$.
\textbf{B}: The travel time $T_{tr}$ decreases with respect to the radius $\eps$. \textbf{C}: $L_{pe}$ increases linearly with the distance $l$. \textbf{D}: $T_{tr}$ increases with the distance fitted by a cubic power law. The amplitude of the current is set to $I=1$.}
\label{fig:fig8}
\end{figure}

To estimate $L_{pe}$ and $T_{tr}$ for a range of different window radii and distances, we solved eq.~\ref{penel} numerically as shown in Fig.~\ref{fig:fig8} and found that $L_{pe}$ is quasi-linear with $l$ but decays smoothly with $\eps$. Our plot also suggests that the dependence of the travel time $T_{tr}$ follows a power law with $\eps$ and $l$. We thus propose empirical expressions for both $L_{pe}$ and $T_{tr}$ that we fitted to our numerical results (Fig.~\ref{fig:fig8}A-D). In summary:
\begin{res}
Our numerical solution suggests that the penetration length $L_{pe}$ is given by the expression
\beq
L_{pe}(l,\eps) \quad = \quad a l - \frac{\eps^2}{l}\,,
\eeq
for $\eps\ll 1$ is the window radius, $l\geq 2\eps$ is the distance between the centers that should not be too large, with
\beq
a \approx 0.8610
\eeq
We fitted the travel time $T_{tr}$ with the power law
\beq
T_{tr}(l,\eps,I) \quad = \quad b\frac{l^3}{I\eps^2}
\eeq
where $I$ is the amplitude of the current, and
\beq
b \approx 1.7445.
\eeq
\end{res}
The previous analysis can be applied under the same condition as above, except that the second window is absorbing. In that case, the solution of Laplace's equation (see appendix) is given by
\begin{align}
&u(x,y,z) = \nonumber \\
& \frac{\pi\eps I}{4} + \eps I\int_0^\infty e^{-mz}\left(J_1(m\eps) J_0\left(m\sqrt{\left(x+\frac{l}{2}\right)^2 + y^2}\right) - \frac{\sin(m\eps)}{2}J_0\left(m\sqrt{\left(x-\frac{l}{2}\right)^2 + y^2}\right) \right) \frac{dm}{m}\,, \label{eq:uplane_abs}
\end{align}
and the integro-differential equation describing the trajectory by
\beq\label{eq:odeabs}
\begin{pmatrix}
\dot{x}(t) \\
\dot{z}(t)
\end{pmatrix}
= \eps I
\begin{pmatrix}
\int_0^\infty e^{-mz}\left(J_1(m\eps)J_1\left(m\left(\frac{l}{2} + x\right)\right) + \frac{\sin(m\eps)}{2}J_1\left(m\left(\frac{l}{2} - x\right)\right)\right)dm \\
\int_0^\infty e^{-mz}\left(J_1(m\eps)J_0\left(m\left(\frac{l}{2} + x\right)\right) - \frac{\sin(m\eps)}{2}J_0\left(m\left(\frac{l}{2} - x\right)\right)\right)dm \\
\end{pmatrix}\,.
\eeq

\begin{figure}[H]
\centering
\includegraphics[width=0.33\linewidth]{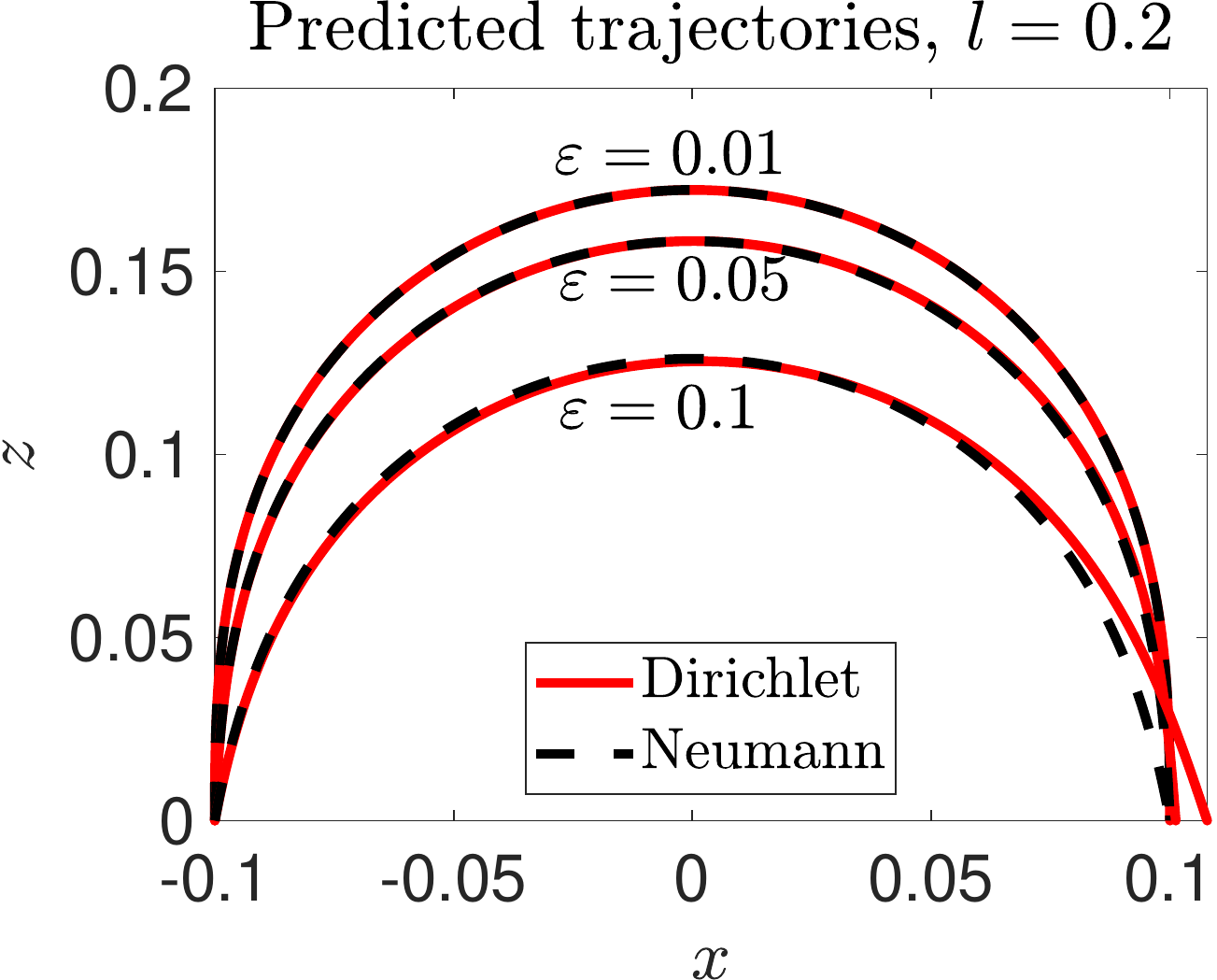}
\caption{\textbf{Exit window with Neumann versus absorbing boundary conditions.} The trajectories are obtained by numerically solving \eqref{eq:ode} and \eqref{eq:odeabs} with the initial condition $(x(0),z(0)) = \left(-\frac{l}{2},0\right)$ using the Matlab ODE solver \textit{ode23} (MathWorks, Natick, MA). The solver is stopped when $z(t)$ reaches $0$ again. Here the amplitude of the current is $I=1$ and the distance is $l=0.2$.}
\label{fig:penetrationLength_abs}
\end{figure}

In Fig.~\ref{fig:penetrationLength_abs} we compare the trajectories obtained by numerically solving eq.~\eqref{eq:ode} and \eqref{eq:odeabs} for $\eps = 0.01,\, 0.05$ and $0.1$. Our results suggest that the choice of boundary conditions does not affect the penetration length, at least for $\eps$ small. However, in contrast to the Neumann case, we observe the end point of the trajectory to slightly drift away from the center of the exiting window when it has absorbing boundary condition. The travel time $T_{tr}$ is also seen to increase.

Finally, for multiple absorbing windows the penetration length can be defined as the maximum of the penetration length between the source and each of the absorbing windows.
\section{Concluding remarks}
The present model and analysis are motivated by the propagation of ions entering through a receptor or a channel and moving inside the cytoplasm of a cell until they reach a pump. In the diffusion approximation, where we modeled pumps or channel as small windows, an influx of ions generates a local difference of concentration between the source and a neighboring target window. Our formula reveals that the concentration difference depends at first order on the size of the window, the current $I$ and the diffusion coefficient $D$, while at second order it also depends on the mean curvature computed at the center of the windows. At third order it depends on the geometrical organization of all the windows via the explicit solution of the Neumann Green's function. A qualitatively similar structure is obtained for the asymptotic expansion of the exiting fluxes. Interestingly at leading-order the exiting flux is the same for all absorbing windows, and here also the geometrical organization has an effect only at third order: the magnitude of the flux is  larger for exiting windows located near the influx receiving window.\\
Then to understand how an influx can perturb the concentration within the domain, we solve the Laplace's equation in the infinite half-space and introduce a novel length scale (the length of penetration) that measures how deep the flow line penetrates inside. The penetration length is independent of the intensity of the field, but depends on the radius and on the distance between the inflow and outflow windows. Our numerical simulations also reveal some scaling laws that should be derived analytically.\\
Finally, this present approach reveals fundamental relations for the biophysics of ionic conduction that elucidate how the concentration could change around voltage-gated channels, a key property for controlling the channel open probability by voltage in dendrites or small protrusions of neuronal cells \cite{holcmanyuste2016}.

\section*{Acknowledgements}
F.P.-L.\ was supported by a postdoctoral fellowship from the Fondation ARC (ARCPDF12020020001505). D.H.\ was supported by the European Research Council (ERC) under the European Union’s Horizon 2020 research and innovation program (grant agreement No 882673).

\begin{appendix}
\section{Appendix: Integral solution of the Laplace's equation with two narrow windows in half-space}\label{sec:canonPDE}
In this appendix we present the integral solutions of Laplace's equation in the three-dimensional half-space that are used in section \S \ref{sec:plength}. Two cases are considered: 1-when an influx and an outflux are imposed on the first and second windows, and 2- when an absorbing boundary condition is imposed on the exiting window.
\subsection{Solution with Neumann boundary conditions on the two windows}
The solution $u(x,y,z)$ of  Laplace's equation
\beq\label{eq:canon}
\frac{\p ^2 u}{\p  x^2} + \frac{\p ^2 u}{\p  y^2} + \frac{\p ^2 u}{\p  z^2} = 0\,, \quad - \infty < x,\,y < \infty\,, \quad z > 0\,,
\eeq
with reflective boundary conditions everywhere on the plane $z=0$ except for a disk of radius $\eps$ on which an influx current is applied, and that is constant at infinity, satisfies the boundary conditions
\beq\label{eq:canonBC1}
\left.-\frac{\p  u}{\p  z}\right|_{z=0} =
\begin{cases}
1 & \sqrt{x^2+y^2} < \eps \\&\\
0 & \sqrt{x^2+y^2} > \eps
\end{cases}\,, \quad \text{and}
\quad u = u_0 \quad \text{for} \quad \sqrt{x^2+y^2+z^2} \to \infty\,, \quad \text{with} \quad z \geq 0\,.
\eeq
The exact solution to \eqref{eq:canon} and \eqref{eq:canonBC1} \cite{crank1975,carslaw1988}, is given by
\beq\label{eq:solcanon1}
u(x,y,z) = u_0 + \eps \int_0^\infty e^{-zm}J_0(m\sqrt{x^2+y^2}) J_1(m\eps) \frac{dm}{m}
\eeq
where $J_n(x)$ are the Bessel functions of order $n$. Using the identity
\beq
\int_0^\infty J_0(mr)J_1(m\eps)dm =
\begin{cases}
0 & r > \eps \\ &\\
\dfrac{1}{\eps} & r < \eps
\end{cases}\,,
\eeq
the solution satisfies the boundary conditions on the $z=0$ plane. The solution decays at infinity like
\beq
u \sim u_0 + \frac{1}{2\rho} + O(\rho^{-3}) \quad \text{as} \quad \rho \to \infty\,,
\eeq
where $\rho = \sqrt{x^2+y^2+z^2}$. For two circular patches of radius $\eps$ whose centers are located at a distance $l$ away on the $y = 0$ line on $\p \Omega$, with the boundary conditions,
\beq
\left.-\frac{\p  u}{\p  z}\right|_{z=0} =
\begin{cases}
1 & \sqrt{\left(x+\frac{l}{2}\right)^2 + y^2} \leq \eps \\&\\
-1 & \sqrt{\left(x-\frac{l}{2}\right)^2 + y^2} \leq \eps \\
0 & \text{elsewhere}
\end{cases}\,,
\eeq
the solution is computed by subtracting the solution of eq.~\eqref{eq:solcanon1} for each window
\beq
u(x,y,z) = \eps \int_0^\infty e^{-mz}J_1(m\eps)\left( J_0\left(m\sqrt{\left(x+\frac{l}{2}\right)^2 + y^2}\right) - J_0\left(m\sqrt{\left(x-\frac{l}{2}\right)^2 + y^2}\right) \right) \frac{dm}{m}\,.
\eeq
\subsection{Solution with an absorbing boundary condition on the exiting window}
When an absorbing boundary condition is imposed on a disk of radius $\eps$ centered at the origin
\beq\label{eq:canonBC2}
\left.\frac{\p  u}{\p  z}\right|_{z=0} = 0\,, \quad \text{for} \quad \sqrt{x^2+y^2} > \eps\,, \quad \text{and} \quad u|_{z = 0} = 0 \quad \text{for} \quad \sqrt{x^2+y^2} < \eps\,,
\eeq
on the plane $z=0$ with $u=u_0$ at infinity, i.e.\ as $\rho = \sqrt{x^2+y^2+z^2} \to \infty$, then $u(x,y,z)$ satisfying eq.~\eqref{eq:canon} is solution of the classical Weber's electrified disk problem (1874) \cite{crank1975,cheviakov2010,gomez2015}. The solution of \eqref{eq:canon} and \eqref{eq:canonBC2} is represented by
\beq\label{eq:solcanon2}
u(x,y,z) = u_0\left(1 - \frac{2}{\pi}\int_0^\infty e^{-zm}J_0(m\sqrt{x^2+y^2}) \sin(m\eps) \frac{dm}{m} \right)\,,
\eeq
and decays at infinity like
\beq
u = u_0\left(1 - \frac{2}{\pi\rho} + O(\rho^{-3}) \right)\,, \quad \text{as} \quad \rho \to \infty\,.
\eeq
Using the identity
\beq
\int_0^\infty J_0(m\sqrt{x^2+y^2}) \sin(m\eps) \frac{dm}{m} = \frac{\pi}{2}\,, \quad \text{for} \quad \sqrt{x^2 + y^2} < \eps\,,
\eeq
the solution \eqref{eq:solcanon2} satisfies $u=0$ on the disk $\p  \Omega_{\eps} = \{(x,y,z)\, |\, \sqrt{x^2+y^2} < \eps \,,\, z = 0\}$. The flux across the boundary of the window is
\beq
\left.\frac{\p  u}{\p  z}\right|_{z=0} = \frac{2u_0}{\pi}\int_0^\infty J_0(m\sqrt{x^2+y^2}) \sin(m\eps) dm =
\begin{cases}
\dfrac{2u_0}{\pi} \dfrac{1}{\sqrt{\eps^2 - \sqrt{x^2+y^2}}} & \sqrt{x^2+y^2} < \eps \\
0 & \sqrt{x^2+y^2} > \eps
\end{cases}\,,
\eeq
By extending the boundary conditions to
\beq
\left.-\frac{\p  u}{\p  z}\right|_{z=0} =
\begin{cases}
1 & \sqrt{\left(x+\frac{l}{2}\right)^2 + y^2} \leq \eps \\
0 & \sqrt{\left(x+\frac{l}{2}\right)^2 + y^2} > \eps \quad \text{and} \quad \sqrt{\left(x-\frac{l}{2}\right)^2 + y^2} > \eps
\end{cases}\,,
\eeq
and
\beq
u = 0 \,, \quad \text{for} \quad \sqrt{\left(x-\frac{l}{2}\right)^2 + y^2} \leq \eps\,,
\eeq
we can add solutions \eqref{eq:solcanon1} to \eqref{eq:solcanon2}, and we get
\begin{align}
u(x,y,z) &= \eps \int_0^\infty e^{-zm}J_0\left((m\sqrt{\left(x+\frac{l}{2}\right)^2+y^2}\right) J_1(m\eps) \frac{dm}{m} \nonumber \\
& + u_0\left(1 - \frac{2}{\pi}\int_0^\infty e^{-zm}J_0\left(m\sqrt{\left(x-\frac{l}{2}\right)^2+y^2}\right) \sin(m\eps) \frac{dm}{m} \right)\,. \label{eq:uplane_abs_temp}
\end{align}
Here the unknown constant $u_0$ is evaluated from the constraint that the total outflux must cancel the total influx, thus giving
\beq
0 = -\frac{2u_0}{\pi}2\pi \int_0^\eps \frac{rdr}{\sqrt{\eps^2 - r^2}} + \pi \eps^2 = -4u_0 \eps + \pi\eps^2\,,
\eeq
from which we readily obtain $u_0 = \frac{\pi\eps}{4}$. Therefore, \eqref{eq:uplane_abs_temp} becomes
\begin{align}
&u(x,y,z) = \nonumber \\
& \frac{\pi\eps}{4} + \eps \int_0^\infty e^{-mz}\left(J_1(m\eps) J_0\left(m\sqrt{\left(x+\frac{l}{2}\right)^2 + y^2}\right) - \frac{\sin(m\eps)}{2}J_0\left(m\sqrt{\left(x-\frac{l}{2}\right)^2 + y^2}\right) \right) \frac{dm}{m}\,,
\end{align}
which we used in eq.~\eqref{eq:uplane_abs}.
\end{appendix}


\normalem
\bibliographystyle{ieeetr}
\bibliography{ref}
\end{document}